\documentclass{article}
\usepackage{amsmath, amssymb, graphics}

\usepackage{amsthm}
\newtheorem{thm}{Theorem}[section]
\newtheorem{df}{Definition}[section]
\newtheorem{rem}{Remark}[section]
\newtheorem{ex}{Example}[section]
\newtheorem{prp}{Proposition}[section]
\newtheorem{lem}{Lemma}[section]
\newtheorem{cor}{Corollary}[section]

\newcommand{\mathsym}[1]{{}}
\newcommand{\unicode}[1]{{}}

\makeatletter
    
    \@addtoreset{equation}{section}
  \makeatother
  
\begin{document}
\title{On the moduli space of polygons with area center}
\author{Fumio HAZAMA\\Tokyo Denki University\\Hatoyama, Hiki-Gun, Saitama JAPAN\\
e-mail address:hazama@mail.dendai.ac.jp\\Phone number: (81)49-296-2911}
\date{\today}
\maketitle
\thispagestyle{empty}

\begin{abstract}
A point $p$ is said to be an area center of a polygon if all of the triangles composed of $p$ and its edges have one and the same area. We construct a moduli space $AC_n$ of such $n$-gons and study its geometry and arithmetic. For every $n\geq 5$, the moduli space is proved to be a rational complete intersection subvariety in $\mathbb{A}^n$. With the help of some subvarieties of low degree in $AC_n$, we also find a unified method of construction of good-looking polygons with area center.\\
keywords: polygon; area center; Chebyshev variety; rational point
\end{abstract}

\section{Introduction}
Let $g$ be a barycenter of a triangle. Then the areas of three triangles made up of $g$ and its edges are one and the same, as is shown by an elementary argument. In view of this fact, $g$ is entitled to be called the area center of the triangle. A general $n$-gon, however, does not always have an area center (see Proposition 1.3). Hence it will be natural to consider when an $n$-gon has an area center. This seemingly innocent problem is, unexpectedly, found to have a intimate connection with the theory of Chebyshev varieties $V_n$ developed in [2], [3]. In the present paper, we investigate the geometry and arithmetic of the moduli space, called $AC_n$, of $n$-gons with the origin as an area center. Among other things we show that $AC_n$ is a rational complete intersection variety of codimension three in $\mathbb{A}^n$. This fact is proved by constructing a Groebner basis for its defining ideal, which in turn provides us with a parametrization of simple form for any $n$. Every parameter, however, does not always correspond to an $n$-gon with an area center of good shape. Here again we find that the family of linear subvarieties of a Chebyshev variety, constructed in [2], plays a crucial role to specify good-looking $n$-gons. Actually we need a slight generalization of the construction, and we introduce the set of strings of parentheses as well as those with bra-ket introduced by Dirac. The latter has no direct connection with quantum theory in this paper, but it will give us a unified viewpoint to study the set of subvarieties of low degree of $AC_n$.\\

The plan of this paper is as follows. Section one gives a precise definition of $n$-gons with an area center, and shows how our problem is related to the theory of Chebyshev varieties. Several useful identities are recalled and generalized, and help us to derive the defining equation of the moduli space $AC_n$ of those $n$-gons. In Section two we construct a Groebner basis of the defining ideal of $AC_n$ and show that it is a rational complete intersection variety. Furthermore we show that $AC_n$ is nonsingular when $n$ is not divisible by four, and in the latter case it has the origin in $\mathbb{A}^n$ as the unique singular point. From Section three on we introduce sets of strings of round brackets including an angle bracket, a bra-ket, a triple bra-ket, and a quadruple bra-ket. In Section three we recall and refine our construction of linear subvarieties of the Chebyshev varieties in [3]. We introduce a partial order on the set of strings of brackets and define several fundamental notions. The ideas of content and of associate polynomials of a string, introduced here, will play a central role throughout the paper. Section four deals with the strings with angle brackets, and Section five deals with those with a bra-ket. The former gives subvarieties of $V_{2n-1}$, and the latter those of $V_{2n}$. Section six together with Section seven is devoted to the construction of subvarieties of $AC_n$ with $n\equiv 0\pmod 4$. Since $AC_n$ is defined to be the intersection of three Chebyshev varieties, we need to make two other strings from a given string. The first is made by an associative transformation (Section six) and the second by a bra-ketting transformation (Section seven). Section eight together with Section nine is devoted to the construction of subvarieties of $AC_n$ with $n\equiv 1\pmod 2$. In these sections we introduce strings with a triple bra-ket and  the associative transformation of the second kind. Furthermore Secion ten together with Section eleven deals with the construction of subvarieties of $AC_n$ with $n\equiv 2\pmod 4$. Here we introduce strings with a quadraple bra-ket and the associative transformation of the third kind. In the final section we conclude the paper by finding that the most symmetric polygons are associated to the set of invariant elements under the action of the symmetric group on $AC_n$.

\section{Polygon with area center}
Let a polygon $P$ in $\mathbb{R}^2$ have $n$ vertices $p_i$ for $i=0,\cdots,n-1$, where $p_k$ is understood to be equal to $p_{k\bmod n}$ when $k$ is smaller than $0$ or greater than $n-1$. We introduce the notion of area-center of $P$ as follows:

\begin{df}
A point $c\in\mathbb{R}^2$ is said to be an area center of $P$ if all the areas of the $n$ triangles with vertices $c,p_i,p_{i+1}$ for $i=0,\cdots,n-1$ are one and the same nonzero real number.
\end{df}

\noindent
\begin{rem}
\normalfont
As is the case in this proposition and throughout the paper, the area means the {\it signed} area.
\end{rem}

\noindent
When $P$ is a triangle, the barycenter $q$ of $P$ coincides with its area-center, since the length of perpendicular from $q$ to each side $p_ip_{i+1}\hspace{1mm}(0\leq i\leq 2)$ is one third of that of perpendicular from $p_{i-1}$. When $n\geq 4$, however, there need not exist an area-center of a general $n$-gon, as is seen later. The main purpose of this paper is to investigate when an $n$-gon has an area-center. In order to begin our study, we express the condition that the origin is an area center of $P$ in terms of the determinant $[p_i,p_{i+1}]\hspace{1mm}(0\leq i\leq n-1)$ of 2-by-2 matrix $(p_i,p_{i+1})$ composed of $p_i$ and $p_{i+1}$:

\begin{prp}
The origin of $\mathbb{R}^2$ is an area-center of $P$ if and only if the equality 
\begin{eqnarray}
[p_{i-1},p_i]=[p_i,p_{i+1}]
\end{eqnarray}
holds for any $i$ with $1\leq i\leq n$.
\end{prp}

\noindent
{\it Proof.} This follows directly from the definition of the area-center, since $[p_{i-1},p_i]$ is twice the area of the triangle composed of $p_{i-1}$, the origin, and $p_i$. \qed\\

\noindent
Since we can express the difference $[p_{i-1},p_i]-[p_i,p_{i+1}]$ as $[p_{i-1}+p_{i+1},p_i]$, the condition (1.1) says that $p_{i-1}+p_{i+1}$ and $p_i$ are linearly dependent. Since $p_i$ could not coincide with the area center by the definition, the condition  (1.1) says that there exists a constant $a_i$ such that 
\begin{eqnarray}
p_{i-1}+p_{i+1}=a_ip_i
\end{eqnarray}
for every $i$. Here we note that the equality (1.2) is expressed as
\begin{eqnarray}
\begin{pmatrix}
p_{i-1} & p_i
\end{pmatrix}
\begin{pmatrix}
0 & -1\\
1 & a_i
\end{pmatrix}
=
\begin{pmatrix}
p_i& p_{i+1}
\end{pmatrix}.
\end{eqnarray}
Hence if we denote the matrix 
$\begin{pmatrix}
0 & -1\\
1 & a_i
\end{pmatrix}$by $A(a_i)$, we see that the equality
\begin{eqnarray}
\begin{pmatrix}
p_0 & p_1
\end{pmatrix}
A(a_1)A(a_2)\cdots A(a_n)=
\begin{pmatrix}
p_0& p_1
\end{pmatrix}.
\end{eqnarray}
holds. Since $p_0$ and $p_1$ are linearly independent by the definition, the equality (1.4) implies that
\begin{eqnarray}
A(a_1)A(a_2)\cdots A(a_n)=
\begin{pmatrix}
1 & 0\\
0 & 1
\end{pmatrix}.
\end{eqnarray}
Conversely if the condition (1.5) is met, then the points $p_0,\cdots,p_{n-1}$ consitute the vertices of an $n$-gon with the origin as an area-center. Furthermore it is an unexpected coincidence that the product on the left hand side of (1.5) plays an important role in our previous study of the Chebyshev varieties in [], [], []. Let us recall some of definitions introduced there. For $n$ independent variables $x_1,\cdots,x_n$, let
\begin{eqnarray*}
U(x_1,\cdots,x_n)=
\begin{pmatrix}
x_1 & -1 & 0 & \cdots & 0 & 0\\
-1 & x_2 & -1 & \cdots & 0 & 0\\
0 & -1 & x_3 & \ddots & 0 & 0\\
\vdots & \vdots & \ddots & \ddots & \ddots & \vdots \\
0 & 0 & 0 & \ddots & x_{n-1} & -1\\
0 & 0 & 0 & \cdots & -1 & x_n
\end{pmatrix},
\end{eqnarray*}
and let 
\begin{eqnarray*}
u(x_1,\cdots, x_n)=\det U(x_1,\cdots,x_n).
\end{eqnarray*}
The zero locus $V(u)$ of $u$ in $\mathbb{A}^n$ is called the {\it Chebyshev variety} (of the second kind), and is denoted simply by $V_n$. For simplicity we use the notation $u[i,j]$ to express $u(x_i,x_{i+1},\cdots,x_j)$ when $i<j$, with the convention that
\begin{eqnarray*}
u[i,i]=x_i,u[i,i-1]=1,
\end{eqnarray*}
for any $i\geq 1$. We recall some identities proved in [], which will be used throughout the paper. The first one is a recurrence relation
\begin{eqnarray}
u[1,n]=x_1u[2,n]-u[3,n],
\end{eqnarray}
which holds for any integer $n\geq 2$. By symmetry, for any $n\geq 2$, we have
\begin{eqnarray}
u[1,n]=x_nu[1,n-1]-u[1,n-2].
\end{eqnarray}
The following identity, which holds for any variables $x,y$, will play a crucial role when we employ an inductive argument in this paper:
\begin{eqnarray}
A(x)A(0)A(y)=-A(x+y).
\end{eqnarray}
Of fundamental importance for our study is the identity:
\begin{eqnarray}
A(x_n)A(x_{n-1})\cdots A(x_1)=
\begin{pmatrix}
-u[2,n-1] & -u[1,n-1]\\
u[2,n] & u[1,n]
\end{pmatrix}.
\end{eqnarray}
Comparing the order of multiplications on the left hand side of (1.9) with that of (1.5), we need to know what occurs if we reverse the order of multiplications on the left hand side:

\begin{lem}
\begin{eqnarray*}
A(x_1)A(x_2)\cdots A(x_n)=
\begin{pmatrix}
-1 & 0 \\
0 & 1
\end{pmatrix}
{}^t(A(x_n)A(x_{n-1})\cdots A(x_1))
\begin{pmatrix}
-1 & 0 \\
0 & 1
\end{pmatrix}.
\end{eqnarray*}
\end{lem}

\noindent
{\it Proof}. This is a direct consequence of the identity 
\begin{eqnarray*}
{}^tA(x)=
\begin{pmatrix}
-1 & 0 \\
0 & 1
\end{pmatrix}
A(x)
\begin{pmatrix}
-1 & 0 \\
0 & 1
\end{pmatrix}. \qed\\
\end{eqnarray*}
It follows from this lemma and the equality (1.9) that
\begin{eqnarray}
A(x_1)A(x_{2})\cdots A(x_n)=
\begin{pmatrix}
-u[2,n-1] & -u[2,n]\\
u[1,n-1] & u[1,n]
\end{pmatrix}.
\end{eqnarray}
Here we record the following polynomial identity which follows from (1.8) and (1.10):

\begin{lem}
For any positive integer $n$ for any $k$ with $2\leq k\leq n-1$, we have
\begin{eqnarray*}
&&u(x_1,\cdots,x_{k-1},0,x_{k+1},\cdots, x_n)\\
&&=-u(x_1,\cdots,x_{k-2},x_{k-1}+x_{k+1},x_{k+2},\cdots, x_n).
\end{eqnarray*}
\end{lem}

\noindent
Now combining (1.5) with (1.10), we obtain the following:

\begin{prp}
An $n$-gon in $\mathbb{R}^2$ has the origin as its area-center if and only if there exists a real solution of the simultaneous equation
\begin{eqnarray}
&&u[1,n]=1, \\
&&u[1,n-1]=0, \\
&&u[2,n]=0.
\end{eqnarray}
\end{prp}

\begin{rem}
\normalfont
Note that the first equality (1.11) can be replaced by the equality 
\begin{eqnarray}
u[2,n-1]=-1, 
\end{eqnarray}
since $\det A(x)=1$ as a polynomial in $x$. 
\end{rem}

\noindent
\begin{df}
We call the subvariety of $\mathbb{A}^n$ defined by the three equations in Proposition {\rm 1.2} the AC-variety and denote it by $AC_n$.
\end{df}

\noindent
We illustrate Proposition 1.2 by a few examples.

\begin{ex}
\normalfont
When $n=3$, the three equations (1.11), (1.12), (1.13) becomes
\begin{eqnarray*}
\left\{
\begin{array}{l}
u[1,3]=x_1x_2x_3-x_1-x_3=1, \\
u[1,2]=x_1x_2-1=0, \\
u[2,3]=x_2x_3-1=0.
\end{array}
\right.
\end{eqnarray*}
Inserting $x_1x_2=1$, which is implied by the second equation, into the first equation, we have $x_1=-1$, and it follows that $x_2=x_3=-1$. Hence the equalities (1.2) for $i=1,2,3$ require one and the same equation $p_0+p_1+p_2=0$, which means that the origin is the barycenter, as is seen before.
\end{ex}

\begin{ex}
\normalfont
When $n=4$, recalling that the equation (1.11) can be replaced by (1.14) (see Remark 1.2), we see that the defining equations of $AC_4$ become
\begin{eqnarray}
\left\{
\begin{array}{l}
u[2,3]=x_2x_3-1=-1, \\
u[1,3]=x_1x_2x_3-x_1-x_3=0, \\
u[2,4]=x_2x_3x_4-x_2-x_4=0.
\end{array}
\right.
\end{eqnarray}
The first equality in (1.16) implies that $x_2$ or $x_3$ is equal to zero. When $x_2=0$, it follows from the second and the third equalities that $x_1+x_3=x_4=0$, and when $x_3=0$, we have $x_1=x_2+x_4=0$. In the case of $x_2=0$, therefore, the four conditions for $i=1,2,3,4$ of (1.2) become
\begin{eqnarray}
\left\{
\begin{array}{l}
p_0+p_2=x_1p_1, \\
p_1+p_3=0, \\
p_2+p_0=-x_1p_3, \\
p_3+p_1=0.
\end{array}
\right.
\end{eqnarray}
Note that the third (resp. the fourth) equation is equivalent to the first (resp. the second) equation, hence (1.16) is simplified to 
\begin{eqnarray}
\left\{
\begin{array}{l}
p_0+p_2=x_1p_1, \\
p_1+p_3=0.
\end{array}
\right.
\end{eqnarray}
Similarly in the case of $x_3=0$, the conditions are reduced to
\begin{eqnarray}
\left\{
\begin{array}{l}
p_0+p_2=0, \\
p_1+p_3=x_2p_2, \\
\end{array}
\right.
\end{eqnarray}
Hence we see that a quadrilateral $P=(p_0,p_1,p_2,p_3)$ has an area-center (not necessarily equal to the origin) if and only if 
\begin{eqnarray}
\mbox{the midpoint of the segment $p_0p_2$ lies on the line $p_1p_3$,} 
\end{eqnarray}
or 
\begin{eqnarray}
\mbox{the midpoint of the segment $p_1p_3$ lies on the line $p_0p_2$. }
\end{eqnarray}
\end{ex}

Using these criteria (1.19), (1.20), we can show that there is a quadrilateral {\it without} area-center:

\begin{prp}
Let $p_0=(0,1), p_1=(-1,0), p_2=(0,-1)\in\mathbb{R}^2$. For a point $p_3\in\mathbb{R}^2$ to constitute a quadrilateral $P=(p_0,p_1,p_2,p_3)$ with area-center, it is necessary and sufficient that $p_3\in V(x-1)\cup (V(y)\setminus\{p_1\}))$ .
\end{prp}

\noindent
{\it Proof of Proposition} 1.3. If the condition (1.19) is satisfied, then the origin must lies on the line $p_1p_3$, hence the point $p_3$ must lie on $x$-axis $V(y)$ and cannot coincides with $p_1$ by Definition 1.1. If the condition (1.20) is satisfied, then the point $p_3$ must lie on the line $V(x-1)$. \qed\\\\
It follows that an $n$-gon with $n\geq 4$ does not always have an area-center. This fact motivates us to investigate what kind of polygons have area-centers.

\section{Geometry of the AC-varieties}
Since the equations (1.11)-(1.13) of $AC_n$ can be defined over any field, we will fix an arbitrary field $K$ and regard $AC_n$ is defined over $K$ from now on. In this section we focus on some of algebro-geometric properties of the AC-varieties. \\

Firstly we construct a Groebner basis of its defining ideal:

\begin{prp}
For any integer $n\geq 4$, let $I=\langle u[1,n]-1, u[2,n], u[1,n-1]\rangle$, the defining ideal of $AC_n$. With respect to the lexicographic order with $x_1>x_2>\cdots>x_n$, a Groebner basis of $I$ is given by
\begin{eqnarray}
I=\langle u[3,n]+1, u[4,n]+x_2, u[3,n-1]+x_1\rangle.
\end{eqnarray}
\end{prp}

\noindent
{\it Proof}. Let $J$ denote the ideal on the right hand side of (2.6). First we show that $I=J$.\\
$J\subset I$: As for the first generator $u[3,n]+1$ of $J$, since $u[1,n]=x_1u[2,n]-u[3,n]$ holds by (1.6), we have
\begin{eqnarray*}
u[3,n]+1=x_1u[2,n]-(u[1,n]-1)\in I
\end{eqnarray*}
This implies in turn that
\begin{eqnarray*}
u[2,n]=x_2u[3,n]-u[4,n]\equiv -x_2-u[4,n]\pmod I,
\end{eqnarray*}
hence the second generator $u[4,n]+x_2$ of $J$ belongs to $I$. As for the third generator, we employ the equality (1.10). Taking the determinants of both sides of (1.10) and noting that $\det A(x)=1$, we have
\begin{eqnarray}
-u[1,n]u[2,n-1]+u[1,n-1]u[2,n]=1.
\end{eqnarray}
Since $u[1,n]\equiv 1 \pmod I, u[2,n]\equiv 0\pmod I$, it follows that 
\begin{eqnarray*}
-u[2,n-1]\equiv 1 \pmod I,
\end{eqnarray*}
which implies that 
\begin{eqnarray*}
u[2,n-1]+1\in I.
\end{eqnarray*}
Hence we have $u[1,n-1]=x_1u[2,n-1]-u[3,n-1]\equiv -x_1-u[3,n-1]\pmod I$. Since $u[1,n-1]\in I$, it follows that $u[3,n-1]+x_1\in I$. Thus we see that $J\subset I$.\\
$I\subset J$: For the second generator $u[2,n]$ of $I$, the equality
\begin{eqnarray*}
u[2,n]=x_2u[3,n]-u[4,n]=x_2(u[3,n]+1)-(u[4,n]+x_2)
\end{eqnarray*}
implies that $u[2,n]\in J$. Therefore we have $u[1,n]-1=x_1u[2,n]-(u[3,n]+1)\in J$. Then the equality (2.7) modulo $J$ this time implies that
\begin{eqnarray*}
-u[2,n-1]\equiv 1\pmod J.
\end{eqnarray*}
hence we have $u[2,n-1]+1\in J$. It follows that
\begin{eqnarray*}
u[1,n-1]=x_1(u[2,n-1]+1)-(u[3,n-1]+x_1)\in J.
\end{eqnarray*}
Thus all of the generators of $I$ belong to $J$, and we see that $I=J$. What remains is to show that the three generators actually gives us a Groebner basis. Note that the leading monomials of them are 
\begin{eqnarray*}
{\rm LM}(u[3,n]+1)&=&x_3x_4\cdots x_n,\\
{\rm LM}(u[4,n]+x_2)&=&x_2,\\
{\rm LM}(u[3,n-1]+x_1)&=&x_1.
\end{eqnarray*}
Since these are pairwise relatively prime, it follows from Buchberger's simplified criterion (see [1, Chapter 2, Section 9, Theorem 3, and Proposition 4], for example) that the three polynomials on the right hand side of (2.1) constitute a Groebner basis of $I$. This completes the proof of Proposition 2.1. \qed\\

\noindent
As a direct consequence, we have the following:
\begin{cor}
The AC-variety $AC_n$ is a complete intersection of dimension $n-3$ in $\mathbb{A}^n$.
\end{cor}

In order to find when $AC_N$ is nonsingular, we first investige the singular points on the hypersurface $V(u[1,n]+1)\subset \mathbb{A}^n$ related to the first element of the Groebner basis (2.1) of the defining ideal. We need the following simple lemma:

\begin{lem}
For any integer $k$ with $1\leq k\leq n$, we have
\begin{eqnarray*}
\frac{\partial}{\partial x_k}u[1,n]=u[1,k-1]u[k+1,n].
\end{eqnarray*}
\end{lem}

\noindent
{\it Proof}. This is a direct consequence of the definition $u[1,n]=\det U(x_1,\cdots,x_n)$ and the rule for the differentiation of a deteminant. \qed\\

\noindent

\begin{prp}
{\rm (1)} If $n\not\equiv 2\pmod 4$, then $V(u[1,n]+1)$ is nonsingular.\\
{\rm (2)} If $n\equiv 2 \pmod 4$, then $V(u[1,n]+1)$ has the origin as its unique singular point.
\end{prp}

\noindent
{\it Proof}. It follows from Lemma 2.1 that a singular point of $V(u[1,n]+1)$ satisfies the simultaneous equation
\begin{eqnarray}
\left\{
\begin{array}{l}
u[2,n]=0,\\
u[1,1]u[3,n] = 0,\\
u[1,2]u[4,n] = 0,\\
\cdots ,\\
u[1,n-2]u[n,n] = 0,\\
u[1,n-1] = 0
\end{array}
\right.
\end{eqnarray}
Here we employ the following simple lemma, which facilitates our proof greatly:
 
\begin{lem}
{\rm (1)} For any $k$ with $2\leq k\leq n$, we have $V(u[1,k])\cap V(u[1,k-1])=\phi$. \\
{\rm (2)}  For any $k$ with $2\leq k\leq n$, we have $V(u[k,n])\cap V(u[k-1,n])=\phi$. \\
\end{lem}

\noindent
{\it Proof of Lemma}. (1) Contrary to the conclusion, suppose that there exists a point $(x_i)\in \mathbb{A}^n$ such that 
\begin{eqnarray}
(x_l)\in V(u[1,k])\cap V(u[1,k-1]).
\end{eqnarray}
 Since $u[1,k]=x_ku[1,k-1]-u[1,k-2]$ holds by (1.7), the assumption (2.4) implies that $(x_i)\in V(u[1,k-2])$, if $k\geq 3$. This leads us eventually to the conclusion
 \begin{eqnarray*}
 (x_i)\in V(u[1,2])\cap V(u[1,1])
 \end{eqnarray*}
 Recall, however, that $u[1,2]=x_1x_2-1, u[1,1]=x_1$, hence
 \begin{eqnarray*}
 V(u[1,2])\cap V(u[1,1])=\phi.
 \end{eqnarray*}
This contradiction completes the proof of (1). The claim (2) can be proved similarly, by appealing to the fact that $V(u[n-1,n])\cap V(u[n,n])=\phi$. \qed\\
 
\noindent
Now we resume the proof of Proposition 2.3. It follows from this lemma that the first equation of (2.3) forces us to choose the first alternative $u[1,1]=0$ of the consequence "$u[1,1]=0$ or $u[3,n]=0$" of the second equation of (2.3). Which in turn forces us to choose the alternative $u[4,n]=0$ of the third equation of (2.3). Repeating similarly, we find that 
\begin{eqnarray}
\mbox{the number of equation }n\mbox{ in (2.3) must be necessarily even.}
\end{eqnarray}
From here on, we assume that the condition (2.5) is met. Then it follows from (2.3) and Lemma 2.2 that
\begin{eqnarray}
u[2,n]&=&u[4,n]=\cdots =u[n,n]=0,\\
u[1,1]&=&u[1,3]=\cdots = u[1,n-1]=0.
\end{eqnarray}
Since we have $u[2,n]=x_2u[3,n]-u[4,n]$ and $u[3,n]\neq 0$ by  our lemma, the first equality in (2.6) implies that $x_2=0$. Repeating this argument with the remaining equalities in (2.6), we must have
\begin{eqnarray*}
x_2=x_4=\cdots = x_n=0.
\end{eqnarray*}
Similaryly it follows from the equalities (2.7) that
\begin{eqnarray*}
x_1=x_3=\cdots = x_{n-1}=0.
\end{eqnarray*}
Thus we are reduced to checking when the origin of $\mathbb{A}^n$ with $n$ even belongs to $V(u[1,n]+1)$. Since $A(0)=\begin{pmatrix}
0 & -1\\
1 & 0
\end{pmatrix}
$ and $A(0)^2=\begin{pmatrix}
-1 & 0\\
0 & -1
\end{pmatrix}
$,  the (2,2)-entry $u[1,n]$ of the matrix on the right hand side of (1.10) takes the value $-1$ at the origin if and only if $n\equiv 2\pmod 4$. This completes the proof of Proposition 2.3. \qed\\

\noindent
Now we can show the following:

\begin{thm}
For any integer $n\not\equiv 0\pmod 4$ , the AC-variety $AC_n$ is a nonsingular complete intersection. When $n\equiv 0\pmod 4$, the AC-variety $AC_n$ has the origin as its unique singular point.  
\end{thm}

\noindent
{\it Proof}. Let 
\begin{eqnarray}
f_1&=&x_1+u[3,n-1],\\
f_2&=&x_2+u[4,n],\\
f_3&=&u[3,n]+1,
\end{eqnarray}
which constitute a Groebner basis of the defining ideal of $AC_n$ by Proposition 2.1. Furthermore let $M=(m_{ij})_{1\leq i\leq3, 1\leq j\leq n}$ denote the 3-by-$n$ matrix with $m_{ij}=\partial f_i/\partial x_j$. Then we see that the 3-by-2 submatrix $(m_{ij})_{1\leq i\leq3, 1\leq j\leq 2}$ of $M$ equals to
\begin{eqnarray*}
\begin{pmatrix}
1 & 0\\
0 & 1\\
0 & 0
\end{pmatrix}.
\end{eqnarray*}
Hence the set of singular points on $AC_n$ is contained in that of the zero locus $V(f_3)$. We know, however, by Proposition 2.3 that the latter is empty when $n\not\equiv 0\pmod 4$. (Note that $f_3=u[3,n]+1$ and the proposition looks at its zero locus in $\mathbb{A}^{n-2}$.) Hence the first assertion of the theorem follows. Thus we are reduced to considering the case when $n\equiv 0\pmod 4$. in this case, Proposition 2.3 tells us that $m_{3j}=0\hspace{1mm}(1\leq j\leq n)$ only when $(x_3,\cdots, x_n)=(0,\cdots, 0)$. Furthermore if  $(x_3,\cdots, x_n)=(0,\cdots, 0)$, then it follows from (1.6) that $u[4,n]=-u[6,n]=\cdots = u[n,n]=x_n=0$, and $u[3,n-1]=0$ for a similar reason. Hence (2.8) and (2.9) implies that if $(x_i)\in AC_n$ and $(x_3,\cdots, x_n)=(0,\cdots, 0)$, then $x_1=x_2=0$ too. This completes the proof. \qed\\

\noindent
Next we show that the AC-varieties are rational:

\begin{prp}
For any integer $n\geq 5$, the AC-variety $AC_n$ is rational. 
\end{prp}

\begin{rem}
The case when $n=4$ is considered in Example $1.2$ already, and we know that $AC_4$ is the union of two lines $V(x_1, x_3, x_2+x_4)$ and $V(x_2,x_4,x_1+x_3)$ in $\mathbb{A}^4$.
\end{rem}

\noindent
{\it Proof}. Let $x=(x_1,\cdots, x_n)$ be a general point on $AC_n$. Since the equality
\begin{eqnarray*}
u[3,n]=x_3u[4,n]-u[5,n]
\end{eqnarray*}
holds for any $n\geq 5$ and $u[3,n]=-1$ by (2.10), we have
\begin{eqnarray}
x_3=\frac{u[5,n]-1}{u[4,n]}.
\end{eqnarray}
The equality (2.9) implies directly that 
\begin{eqnarray*}
x_2=-u[4,n].
\end{eqnarray*}
Furthermore it follows from (2.8) that
\begin{eqnarray*}
x_1=-u[3,n-1]=-x_3u[4,n-1]+u[5,n-1].
\end{eqnarray*}
Inserting the right hand side of (2.11) into this expression, we have
\begin{eqnarray*}
x_1&=&-\left(\frac{u[5,n]-1}{u[4,n]}\right)u[4,n-1]+u[5,n-1]\\
&=&\frac{-(u[4,n-1]u[5,n]-u[4,n]u[5,n-1])+u[4,n-1]}{u[4,n]}.\\
\end{eqnarray*}
Amazingly enough, the expression in the bracket in the numerator of the last fraction is equal to 
\begin{eqnarray*}
\det(A(x_4)A(x_5)\cdots A(x_n))
\end{eqnarray*}
by (1.10), hence it is equal to one and we have
\begin{eqnarray*}
x_1=\frac{u[4,n-1]-1}{u[4,n]}.
\end{eqnarray*}
This completes the proof of the rationliaty. \qed\\

\section{Brackets}
In [3] we see how a string of brackets corresponds to a subvariety of the Chebyshev varieties. We will generalize the correspondence in various ways in this paper. \\

A string of left brackets and right brackets is said to be balanced if each left bracket has a matching right bracket. To be more precise, we define a function $s:\{(, )\}\rightarrow \{\pm 1\}$ by $s("(")=1, s(")")=-1$. Then the balancedness is defined as follows:

\begin{df}
A string $(b_i)_{1\leq i\leq 2n}$ is balanced if and only if $\sum_{1\leq i\leq 2n}s(b_i)=0$ and 
$\sum_{1\leq i\leq k}s(b_i)\geq 0$ for any $k$ with $1\leq k\leq 2n$.
\end{df}

\noindent
For a positive integer $n$, let $Par_n$ denote the set  of balanced strings of brackets of length $2n$. For any balanced string $b=(b_1,\cdots,b_{2n})\in Par_n$, let  $Pairs(b)$ denote the set of matching pair of left and right brackets. On the set $Pairs(b)$ we introduce a partial order $"\prec"$: For $(b_i,b_j), (b_k,b_{\ell})\in Pairs(b)$, $(b_i,b_j)\prec(b_k,b_{\ell})$ if and only if $k<i$ and $j<\ell$.  We denote by $L(b)$ (resp. $R(b)$) the set of left (resp. right) brackets among $\{b_i\}$. When $b_i$ belongs to $L(b)$, its matching right bracket will be denoted by $r(b_i)\in R(b)$. Furthermore for any $b_i\in L(b)$ let $ls(b_i)$ denote the set of left brackets in $\{b_{i+1},\cdots, r(b_i)\}$. Note that $ls(b_i)=\phi$, when $b_{i+1}\in R(b)$. Association of a left bracket $b_i$ with the matching pair $(b_i,r(b_i))$ defines a natural bijection from the set of left brackets $L(b)$ to the set of matching pairs $Pairs(b)$. Through this correspondence we translate the partial order $"\prec"$ on $Pairs(b)$ into that on $L(b)$, which will be denoted by the same symbol $"\prec"$.   Then the poset $L(b)$ endowed with this partial order becomes a ranked poset with the rank function {\it rank} defined as follows:
\begin{df}
For a balanced bracket $b=(b_1,\cdots,b_{2n})\in Par_n$ and for any left bracket $b_i$ in $b$, the rank  of $b_i$, denoted by $rank(b_i)$, is defined inductively by the rule
\begin{displaymath}
\left\{
\begin{array}{ll}
rank(b_i)=1, & \mbox{if}\hspace{1mm}r(b_i)=b_{i+1},  \\
rank(b_i)=\max_{b_k\in ls(b_i)}rank(b_k)+1, & otherwise.
\end{array}
\right.
\end{displaymath}
Moreover the maximum of the ranks of the left brackets in $L(b)$ is called the {\rm height} of $b$ and is denoted by $height(b)$. Furthermore for any $k\in [1,height(b)]$ we set
\begin{eqnarray*}
ls_{rank=k}(b)=\{b_i\in L(b);rank(b_i)=k\}.
\end{eqnarray*}
\end{df}

\noindent
In order to associate a set of polynomials with a balanced string, we regard a string $b=(b_1,\cdots, b_{2n})\in Par_n$ is situated along the number line such that $b_i$ is located at $i-0.5$ for every $i\in\{1,\cdots,2n\}$. Its coordinate will be denoted by $pos(b_i)$ so that $pos(b_i)=i-0.5$. For any pair of real numbers $p,q$ with $p<q$, we put $[p,q]=\{k\in\mathbb{Z};p\leq k\leq q\}$, and for any $i,j$ with $1\leq i<j\leq 2n$, we put $[b_i,b_j]=[pos(b_i),pos(b_j)]$. Most important notion for a string in $Par_n$ is that of {\it content}. It is defined as follows. Let $b_i$ be a left bracket. Then the content $cont(b_i)$ is defined by the following rule:
\begin{eqnarray*}
cont(b_i)=[b_i,r(b_i)]\setminus \bigcup_{b_k\in ls(b_i)}[b_k,r(b_k)].
\end{eqnarray*}
Furtheremore collecting all contents of $b_i\in L(b)$, we set
\begin{eqnarray*}
cont(b)&=&\bigcup_{b_i\in L(b)}\{cont(b_i)\},\\
num(b)&=&\bigcup_{b_i\in L(b)}cont(b_i).
\end{eqnarray*}
Note here there is a slight distinction between $cont(b)$ and $num(b)$. Let us illustrate the notation introduced above by some examples.\\

\noindent
Example 3.1\\
Let $b=(b_1,\cdots,b_6)="(()())"\in Par_3$. We see that $L(b)=\{b_1,b_2,b_4\}$ and $R(b)=\{b_3,b_5,b_6\}$. The matching right brackets are given by $r(b_1)=b_6, r(b_2)=b_3$, and $r(b_4)=b_5$, hence the set of matching pairs is given by
\begin{eqnarray*}
Pairs(b)=\{(b_1,b_6), (b_2,b_3), (b_4,b_5)\}.
\end{eqnarray*}
The pairs are ordered as
\begin{eqnarray*}
(b_2,b_3)\prec (b_1,b_6)\succ (b_4,b_5),
\end{eqnarray*}
and $(b_2,b_3)$ and $(b_4,b_5)$ are not comparable. The rank of $(b_2,b_3)$ and $(b_4,b_5)$ are equal to 1, and that of $(b_1,b_6)$ is 2. The contents are computed as follows;
\begin{eqnarray*}
cont(b_1)&=&[b_1,b_6]-([b_2,b_3]\cup [b_4,b_5])\\
&=&\{1,2,3,4,5\}-(\{2\}\cup\{4\})\\
&=&\{1,3,5\},\\
cont(b_2)&=&\{2\},\\
cont(b_4)&=&\{4\}.
\end{eqnarray*}
Therefore we have
\begin{eqnarray*}
cont(b)&=&\{\{1,3,5\},\{2\},\{4\}\},\\
num(b)&=&\{1,2,3,4,5\}.
\end{eqnarray*}

\noindent
Example 3.2\\
Let $b=(b_1,\cdots,b_{10})="(())()(())"\in Par_5$. For this string $b$ we have

\begin{eqnarray*}
cont(b)&=&\{\{1,3\},\{7,9\},\{2\},\{5\},\{8\}\},\\
num(b)&=&\{1,2,3,5,7,8,9\}.
\end{eqnarray*}

\noindent
In these examples we observe that each $cont(b_i)$ for $b_i\in L(b)$ consists either wholly of odd integers  or wholly of even integers. The following proposition shows that this observation is correct:

\begin{prp}
For any string $b\in Par_n$ and for any $b_i\in L(b)$, the parities of elements in $cont(b_i)$ are one and the same.
\end{prp}

\noindent
{\it Proof.} We prove this by induction on the rank of $b_i\in L(b)$. When the rank is equal to one, there is nothing to prove since $cont(b_i)=\{i\}$ in this case. Therefore suppose that $rank(b_i)=r\geq 2$ and the assertion is proved for every left bracket of smaller rank. Let us put $ls(b_i)\cap ls_{rank=r-1}(b)=\{b_{i_1},\cdots,b_{i_p}\}$. Note that $i_1=i+1$, since $rank(b_i)=1$ if otherwise. Assume that $i$ is odd. Then we have $\min(cont(b_{i_1}))=\min(cont(b_{i+1}))=i+1\equiv 0\pmod 2$, and hence $\max(cont(b_{i_1}))\equiv 0\pmod 2$ by the induction hypothesis. It follows that the second smallest element of $cont(b_i)$ is equal to $\max(cont(b_{i_1}))+1$, and hence it is odd, which in turn shows that $\min(cont(b_{i_2}))=\max(cont(b_{i_1}))+2\equiv 0\pmod 2$, and hence the third smallest element of $cont(b_i)$ is odd, and so on. Arguing in this way we arrive at the conclusion that every element in $cont(b_i)$ is odd. Since the case when $i$ is even can be treated similarly, we finish the proof. \qed\\

\noindent
Now we attach a set of linear polynomials to each balanced string of brackets. Fix a nonnegative integer $n$ and let $x_1,\cdots x_{2n}$ be independent variables. For any subset $S\subset [1,2n]$, let $f_S$ denote the linear polynomial $\sum_{k\in S}x_k$. When $b\in Par_n$, we define a set of polynomial $f_b$ as
\begin{eqnarray*}
f_b=\{f_{cont(b_i)};b_i\in L(b)\}.
\end{eqnarray*}
Therefore for the string $b$ in Example 1.1, we have
\begin{eqnarray*}
f_b=\{x_1+x_3+x_5,x_2,x_4\},
\end{eqnarray*}
and for $b$ in Example 1.2, we have
\begin{eqnarray*}
f_b=\{x_1+x_3,x_7+x_9,x_2,x_5,x_8\}.
\end{eqnarray*}
\noindent
We will see later that the zero locus $V(f_b)$ for every $b\in Par_n$ provides us with a subvariety of $AC_n$. The following proposition plays a crucial role several times when  we employ an inductive argument to show this fact and some related results: 

\begin{prp}
Let $b=(b_i)\in Par_n$ and suppose that there exists a minimal element $b_k\in L(b)$ which is not maximal. Let $b'\in Par_{n-1}$ denote the string $(b_1,\cdots,b_{k-1},b_{k+2},\cdots,b_{2n})$. Let $p_k:\mathbb{A}^{2n-1}\rightarrow\mathbb{A}^{2n-3}$ denote the linear map defined by the rule
\begin{eqnarray*}
p_k(x_1,\cdots,x_{2n-1})=(x_1,\cdots,x_{k-2},x_{k-1}+x_{k+1},x_{k+2},\cdots, x_{2n-1}).
\end{eqnarray*}
Then we have
\begin{eqnarray}
p_k^*(f_{b'})=f_b\setminus \{x_k\}.
\end{eqnarray}
\end{prp}

\noindent
Before we begin our proof, we illustrate the strings $b$ and $b'$ when $b=(())(()())\in Par_5$ and we choose $b_6$ as a minimal but non-maximal element in $L(b)$. The upper cover of $b_6$ is $b_5$ in this case. The content $cont(b)$ and the associated set of polynomials $f_b$ are given by
\begin{eqnarray*}
cont(b)=\{cont(b_1), cont(b_2), cont(b_5), cont(b_6), cont(b_8)\},
\end{eqnarray*}
with
\begin{eqnarray*}
cont(b_1)&=&\{1,3\},\\ 
cont(b_2)&=&\{2\},\\
cont(b_5)&=&\{5,7,9\},\\
cont(b_6)&=&\{6\}, \\
cont(b_8)&=&\{8\},
\end{eqnarray*}
and
\begin{eqnarray*}
f_b=\{x_1+x_3, x_2, x_5+x_7+x_9, x_6, x_8\}.
\end{eqnarray*}
As for the string $b'$, we have
\begin{eqnarray*}
cont(b')=\{cont(b'_1), cont(b'_2), cont(b'_5), cont(b'_6)\},
\end{eqnarray*}
with
\begin{eqnarray*}
cont(b'_1)&=&\{1,3\},\\ 
cont(b'_2)&=&\{2\},\\
cont(b'_5)&=&\{5,7\},\\
cont(b'_6)&=&\{6\}.
\end{eqnarray*}
and
\begin{eqnarray*}
f_{b'}=\{y_1+y_3, y_2, y_5+y_7, y_6\}.
\end{eqnarray*}
The map $p_6:\mathbb{A}^9\rightarrow\mathbb{A}^7$ is given by
\begin{eqnarray*}
p_6(x_1,\cdots,x_9)=(y_1,\cdots, y_7)
\end{eqnarray*}
with
\begin{eqnarray*}
y_1=x_1, \cdots, y_4=x_4, y_5=x_5+x_7, y_6=x_8, y_7=x_9.
\end{eqnarray*}
Therefore we have
\begin{eqnarray*}
p_6^*(f_{b'})&=&\{x_1+x_3, x_2, x_5+x_7+x_9, x_8\}\\
&=&f_b\setminus\{x_6\},
\end{eqnarray*}
which shows that Proposition 3.2 holds for this $b\in Par_5$. This example will facilitate the reader to follow the argument in the proof.\\

\noindent
{\it Proof}. Let $b_{\ell}\in L(b)$ be the upper cover of $b_k$, which is unique by the definition of the partial order on $L(b)$. It follows from the definition of the content that
\begin{eqnarray*}
cont(b_k)&=&\{k\},\\
cont(b_{\ell})&\ni& k-1, k+1,
\end{eqnarray*}
and that
\begin{eqnarray}
cont(b_m)\not\ni k-1,k, k+1,
\end{eqnarray}
for any other $b_m\in L(b)\setminus \{b_k, b_{\ell}\}$, since the contents for distinct left brackets have no element in common. By the definition of $p_k$, if we employ the coordinates $(x_i)_{1\leq i\leq 2n-1}$ for the source $\mathbb{A}^{2n-1}$ and $(y_j)_{1\leq j\leq 2n-3}$ for the target $\mathbb{A}^{2n-3}$, then $p_k((x_i))=(y_j)$ with
\begin{eqnarray}
y_j=
\left\{
\begin{array}{ll}\nonumber
x_j, & j\leq k-2, \\
x_{k-1}+x_{k+1}, & j=k-1, \\\nonumber
x_{j+2}, & j\geq k.
\end{array}
\right.\\
\end{eqnarray}
We need to know the difference between the content of $b$ and that of $b'$, and will see how the associated sets of polynomials $f_b$ and $f_{b'}$ are related through the pull-back by the linear map $p_k.$ For this pupose we introduce a pair of monotone maps $minus2_k:[1,2n-1]\rightarrow [1,2n-3]$, $plus2_k:[1,2n-3]\rightarrow [1,2n-1]$, defined by the following rules:
\begin{eqnarray*}
minus2_k(i)&=&
\left\{
\begin{array}{ll}
i, & i\leq k-1, \\
k-1, & i=k, \\
i-2, & i\geq k+1,
\end{array}
\right.\\
plus2_k(i)&=&
\left\{
\begin{array}{ll}
i, & i\leq k-1, \\
i+2, & i\geq k.
\end{array}
\right.
\end{eqnarray*}
Note here that the pair satisfies the identity
\begin{eqnarray}
(plus2_k\circ minus2_k)(i)=i,\hspace{1em}for\hspace{1mm}i\neq k,k+1.
\end{eqnarray}
With the help of the map $minus2_k$, we see that the content of $b'$ is expressed as follows:
\begin{eqnarray*}
cont(b')=\{minus2_k(cont(b_j));b_j\in L(b)-\{b_k\}\}.
\end{eqnarray*}
On the other hand if $b_m\in L(b)\setminus \{b_k, b_{\ell}\}$, then $k-1, k, k+1\not\in cont(b_m)$ by (3.2), hence it follows from (3.3) and (3.4) that
\begin{eqnarray*}
p_k^*(f_{minus2_k(cont(b_m))})&=&p_k^*(\sum_{j\in minus2_k(cont(b_m))}y_j)\\
&=&\sum_{j\in minus2_k(cont(b_m))}x_{plus2_k(j)}\\
&=&\sum_{i\in cont(b_m)}x_i\\
&=&f_{b_m}
\end{eqnarray*}
For $b_{\ell}$, recall that $k-1, k+1\in cont(b_{\ell})$, and hence $k, k+2\not\in cont(b_{\ell})$ by Proposition 3.1. Therefore, noting that $minus2_k(k-1)=minus2_k(k+1)=k-1$, we can compute the pull-back as follows:
 \begin{eqnarray*}
p_k^*(f_{minus2_k(cont(b_{\ell}))})&=&p_k^*(\sum_{j\in minus2_k(cont(b_{\ell}))}y_j)\\
&=&p_k^*(\sum_{j\in minus2_k(cont(b_{\ell}))\cap [1,k-2]}y_j)\\
&&+p_k^*(y_{k-1})\\
&&+p_k^*(\sum_{j\in minus2_k(cont(b_{\ell}))\cap [k,2n-3]}y_j)\\
&=&(\sum_{j\in minus2_k(cont(b_{\ell}))\cap [1,k-2]}x_j)\\
&&+(x_{k-1}+x_{k+1})\\
&&+(\sum_{j\in minus2_k(cont(b_{\ell}))\cap [k,2n-3]}x_{j+2})\\
&=&\sum_{j\in cont(b_{\ell})\cap [1,k-2]}x_j\\
&&+(x_{k-1}+x_{k+1})\\
&&+(\sum_{j\in cont(b_{\ell})\cap [k,2n-3]}x_j)\\
&=&f_{b_{\ell}}.
\end{eqnarray*}
Thus we see that the equality (3.1) holds. This completes the proof. \qed\\

\section{Strings with angle brackets}
In this section we introduce a set of strings of length $2n$ each of which consists of $n-1$ pairs of round brackets and a pair of angle brackets $\langle$$\rangle$. Given such a string $b=(b_i)_{1\leq i\leq 2n}$, assume that $b_k=\langle$ and $b_{\ell}=\rangle$. Then the string $b$ is said to be {\it balanced} if all of the three substrings $(b_i)_{1\leq i\leq k-1}, (b_i)_{k+1\leq i\leq \ell-1}$, and $(b_i)_{\ell+1\leq i\leq 2n}$ are balanced in the sense of Definition 3.1. Let $Ang_n$ denote the set of balanced strings consisting of $n-1$ pairs of round brackets and a pair of angle brackets. Accordingly the set of left brackets $L(b)$ of $b\in Ang_n$ consists of $n-1$ left round brackets and one $\langle$ in $b$. A partial order on $L(b)$ is defined in the same way as for $Par_n$. Note that if $b\in Ang_n$ and $b_k=\langle$, then $b_k$ is maximal by the definition of balancedness. For any $b\in Ang_n$. we also define the content $cont(b)$ in the same way as for $Par_n$ with the only exception that when $b_k=\langle$, we express its content by using angle bracket. Let us illustrate this point by the following example.\\\\
\noindent
Example 4.1. Let $b=(())\langle()()\rangle()\in Ang_6$. Its content is given by
\begin{eqnarray*}
cont(b)=\{\{2\},\{6\},\{8\},\{11\}, \{1,3\},\langle 5,7,9\rangle \}.
\end{eqnarray*}\\
Fix an element $c\in K$. For any $b\in Ang_n$,we assign a polynomial to every bracket in $L(b)$ by the following rule: If $b_i$ is an round bracket, then $f_{b_i}$ is defined to be the same as in the previous section, and if $b_i=\langle$, and $cont(b_i)=\langle k_1,\cdots,k_m\rangle$, then we define
\begin{eqnarray*}
f_{b_i}(c)=\sum_{1\leq i\leq m}x_{k_i}-c.
\end{eqnarray*}
Gathering these polynomials, we define the set of polynomials $f_b(c)$ by the rule
\begin{eqnarray}
f_b(c)=\{f_{b_{\ell}};b_{\ell}\in L(b)\setminus \{b_i\}\}\cup \{f_{b_i}(c)\}.
\end{eqnarray}
Therefore for the string $b$ given in Example 4.1, we have
\begin{eqnarray*}
f_b(c)=\{x_2, x_6, x_8, x_{11}, x_1+x_3,x_5+x_7+x_9-c\},
\end{eqnarray*}
The reason why we introduce this kind of the sets of polynomials is that they provide us with a family of linear subvarieties of Chebyshev varieties and their generalizations. For any $c\in K$ and for any pair $(i,j)$ of positive integers with $i<j$, let
\begin{eqnarray*}
V_{[i,j]}(c)=\{(x_i,\cdots, x_j)\in\mathbb{A}^{j-i+1};u[i,j]=c\}.
\end{eqnarray*}

\begin{thm}
For any $c\in K$ and for any $b\in Ang_n$, we have 
\begin{eqnarray*}
V(f_b(c))\subset V_{[1,2n-1]}((-1)^{n-1}c).
\end{eqnarray*}
\end{thm}

\noindent
{\it Proof}.
We prove this by induction on $n$. When $n=1$, the set $Ang_1$ consists solely of $\langle$$\rangle$. If we call this string $b$, then
\begin{eqnarray*}
f_b(c)=\{x_1-c\}.
\end{eqnarray*}
On the other hand, we have 
\begin{eqnarray*}
V_{[1,1]}((-1)^{1-1}c)=\{(x_1)\in\mathbb{A}^1;u[1,1]=c\}=\{(x_1)\in\mathbb{A}^1;x_1=c\},
\end{eqnarray*}
hence the assertion is proved. When $n=2$, the set $Ang_2$ has three strings:
\begin{eqnarray*}
Ang_2=\{()\langle\rangle, \langle\rangle (), \langle ()\rangle\}.
\end{eqnarray*}
Let us call these strings $b^1,b^2,b^3$, respectively. Then it follows from the definition that
\begin{eqnarray*}
f_{b^1}(c)&=&\{x_1,x_3-c\},\\
f_{b^2}(c)&=&\{x_1-c,x_3\},\\
f_{b^3}(c)&=&\{x_2,x_1+x_3-c\}.
\end{eqnarray*}
On the other hand, since $u[1,3]=x_1x_2x_3-x_1-x_3$, we have 
\begin{eqnarray*}
V_{[1,3]}((-1)^{2-1}c)=\{(x_i)\in\mathbb{A}^3;x_1x_2x_3-x_1-x_3=-c\},
\end{eqnarray*}
and we see directly that the assertion holds true. Now assume that $n\geq 3$ and that the assertion is proved for any $n'<n$. Take an arbitrary balanced string $b=(b_1,\cdots,b_{2n})\in Ang_n$. We divide our argument into two cases: (I) when every minimal element in $L(b)$ is maximal, and (II) when there exists an element in $L(b)$ which is minimal but not maximal. In another word, the Hasse diagram for $b$ in the case (I) is totally disconnected and that in the case (II) is not so.\\\\
Case (I) Every minimal element in $L(b)$ is maximal. In this case we have $L(b)=\{b_1,b_3,\cdots,b_{2n-1}\}$ and that there is a unique integer $k$ with $1\leq k\leq n$ such that $b_{2k-1}=\langle$. It follows that
\begin{eqnarray*}
f_b(c)=\{x_1,x_3,\cdots,x_{2k-3},x_{2k-1}-c,x_{2k+1},\cdots,x_{2n-1}\}
\end{eqnarray*}
On the other hand the equalities (1.6) and (1.7) imply that
\begin{eqnarray*}
&&u(0,x_2,0,x_4,\cdots,0,x_{2k-2},x_{2k-1},x_{2k},0,x_{2k+2},\cdots,0, x_{2n-2},0)\\
&=&(-1)^{k-1}\cdot (-1)^{n-k}\cdot u(x_{2k-1})\\
&=&(-1)^{n-1}x_{2k-1}
\end{eqnarray*}
Therefore if $a=(a_1,\cdots,a_{2n-1})\in V(f_b(c))$, then we have
\begin{eqnarray*}
u(a_1,\cdots, a_{2n-1})&=&(-1)^{n-1}a_{2k-1}\\
&=&(-1)^{n-1}c.
\end{eqnarray*}
It follows that $a\in V_{[1,2n-1]}((-1)^{n-1}(c))$. \\\\
Case (II): There exists an element in $L(b)$ which is minimal but not maximal. Therefore we are in the situation treated in Proposition 3.2. As is assumed there, we call such an element $b_k$, and let $b'$ denote the string defined by
\begin{eqnarray*}
b'=(b_1,\cdots,b_{k-1},b_{k+2},\cdots,b_{2n}).
\end{eqnarray*}
Note that $b_k$ is not angle bracket by the definition of balancedness, and hence $b'$ belongs to $Ang_{n-1}$. Let $p_k:\mathbb{A}^{2n-1}\rightarrow\mathbb{A}^{2n-3}$ denote the map defined by
\begin{eqnarray*}
p_k(x_1,\cdots,x_{2n-1})=(x_1,\cdots,x_{k-2},x_{k-1}+x_{k+1},x_{k+2},\cdots,x_{2n-1}),
\end{eqnarray*}
and let $\bar{p}_k$ denote its restriction to $V(x_k)\subset\mathbb{A}^{2n-1}$. Then it follows from Proposition that
\begin{eqnarray*}
p_k^*(f_{b'}(c))=f_b(c)\setminus\{x_k\}.
\end{eqnarray*}
Hence we have
\begin{eqnarray}
p_k^{-1}(V(f_{b'}(c)))=V(f_b(c)).
\end{eqnarray}
Furthermore suppose that $a=(a_i)\in V(x_k)$ and $\bar{p}_k(a)\in V_{[1,2n-3]}((-1)^{n-2}c)$. It follows that 
\begin{eqnarray*}
u(a_1,\cdots,a_{k-2},a_{k-1}+a_{k+1},a_{k+2},\cdots,a_{2n-1})=c.
\end{eqnarray*}
The left hand side is equal to $-u(a_1,\cdots,a_{k-1},0,a_{k+1},\cdots,a_{2k-1})$ by Lemma 1.2. It follows that $u(a)=(-1)^{n-1}c$ and $a\in V_{[1,2n-1]}((-1)^{n-1}c)$. Thus we have
\begin{eqnarray}
p_k^{-1}(V_{[1,2n-3]}((-1)^{n-2}c))&\subset & V(x_k)\cap V_{[1,2n-1]}((-1)^{n-1}c).
\end{eqnarray}
On the other hand, by the induction hypothesis we have 
\begin{eqnarray*}
V(f_{b'}(c))\subset V_{[1,2n-3]}((-1)^{n-2}c),
\end{eqnarray*}
 hence (4.2) and (4.3) imply that
\begin{eqnarray*}
V(f_b(c))\subset V(x_k)\cap V_{[1,2n-1]}((-1)^{n-1}c)\subset V_{[1,2n-1]}((-1)^{n-1}c).
\end{eqnarray*}
This completes the proof for Case (II), hence at the same time the proof of Theorem 4.1. \qed\\\\

Note that when $c=0$, the set $f_b(0)$ coincides with the set $f_b$ introduced in the previous section. Hence, letting $c=0$ in Theorem 4.1, we obtain the following corollary which is proved by a different and complicated method in [3]: 
\begin{cor}
For any $b\in Par_n$, we have $V(f_b)\subset V_{2n-1}$.
\end{cor}

\section{Strings with a bra-ket}
Here {\it bra-ket} means the pair of {\it bra} $\langle |$ and {\it ket} $| \rangle$ as is introduced by Dirac. We regard a bra-ket as consisting of three elements $"\langle"$ (the left angle bracket), $"|"$ (the vertical bar), and $"\rangle"$ (the right angle bracket). We define the balancedness of a string $b=(b_1,\cdots, b_{2n+1})$ with $n-1$ pairs of round brackets and one bra-ket as follows. Let $b_i=\langle, b_j=|, b_k=\rangle$ with $i<j<k$. Then $b$ is said to be balanced if the four remaining substrings$(b_1,\cdots,b_{i-1}), (b_{i+1},\cdots, b_{j-1}), (b_{j+1},\cdots, b_{k-1})$, and $(b_{k+1},\cdots, b_{2n+1})$ are balanced in the sense of Definition 3.1. We employ the convention that the left angle bracket $b_i$ matches with the vertical bar $b_j$ and the same vertical bar $b_j$ matches with the right angle bracket $b_k$. Accordingly we regard a vertical bar both as a left bracket and as a right bracket, and we can define a partial order on the set $L(b)$ of left brackets in $b$ in the same way as for $Par_n$. We denote the set of balanced strings of length $2n+1$ with $n-1$ pairs of round brackets and with one bra-ket by $Bra_n$. The content of a string $b\in Bra_n$ is defined in the same way as for $Par_n$ with the only difference that the content of the bra-ket $(b_i,b_j,b_k)$ is denoted by $\langle cont(b_i) | cont(b_j) \rangle$. The following example should clarify the meaning.\\\\
Example 5.1.\\
Let $b=()\langle()|()()\rangle()\in Bra_6$. The set $L(b)$ (resp. $R(b)$) of left (resp. right) brackets is given by
\begin{eqnarray*}
L(b) &=&\{b_1,b_3,b_4,b_6,b_7,b_9,b_{12}\},\\
R(b)&=&\{b_2,b_5,b_6,b_8,b_{10},b_{11},b_{13}\}.
\end{eqnarray*}
The content of $b$ is given by
\begin{eqnarray*}
cont(b)=\{\{1\},\{4\},\{7\},\{9\},\{12\},\langle 3,5|6,8,10\rangle\}.
\end{eqnarray*}
\\
\noindent
For a string $b=(b_1,\cdots, b_{2n+1})\in Bra_n$, let $b_i=\langle, b_j=|, b_k=\rangle$. For an arbitrary $c\in K$, we define the polynomial $f_{b_j}(c)$ by
\begin{eqnarray}
f_{b_j}(c)=\sum_{p\in cont(b_i)}x_p\sum_{q\in cont(b_j)}x_q-c,
\end{eqnarray}
and for a left round bracket $b_{\ell}\in L(b)$, the polynomial $f_{b_{\ell}}$ is defined to be the same as in Section three. Gathering these polynomials, we define the set of polynomials $f_b(c)$ by the rule
\begin{eqnarray}
f_b(c)=\{f_{b_{\ell}};b_{\ell}\in L(b)\setminus \{b_i,b_j\}\}\cup \{f_{b_j}(c)\}.
\end{eqnarray}
Therefore for the string $b$ in Example 5.1, we have
\begin{eqnarray*}
f_b(c)=\{x_1,x_4,x_7,x_9,x_{12},(x_3+x_5)(x_6+x_8+x_{10})-c\}.
\end{eqnarray*}
The reason why we introduce this third kind of polynomials is that the zero locus of $f_b(c)$ for $b\in Bra_n$ gives rise to a subvariety of degree two of the Chebyshev varieties. More precisely we can prove the following:
\begin{thm}
For any $b\in Bra_n$, we have $V(f_b(c))\subset V_{[1,2n]}((-1)^{n-1}(c-1))$.
\end{thm}
\noindent
{\it Proof}.
We prove this by induction on $n$. When $n=1$, the set $Bra_1$ consists solely of $\langle$$|$$\rangle$. If we call this string $b$, then its content is given by
\begin{eqnarray*}
cont(b)=\{\langle1|2\rangle\}.
\end{eqnarray*}
Therefore we have
\begin{eqnarray*}
f_b(c)=\{x_1x_2-c\}
\end{eqnarray*}
by (5.1) and (5.2). On the other hand, we have 
\begin{eqnarray*}
V_{[1,2]}((-1)^{1-1}(c-1))&=&\{(x_1,x_2)\in\mathbb{A}^2;u[1,2]=c-1\}\\
&=&\{(x_1,x_2)\in\mathbb{A}^2;x_1x_2-1=c-1\}\\
&=&\{(x_1,x_2)\in\mathbb{A}^2;x_1x_2-c=0\},
\end{eqnarray*}
hence the assertion holds. When $n=2$, the set $Bra_2$ has four strings:
\begin{eqnarray*}
Bra_2=\{()\langle |\rangle, \langle |\rangle (), \langle ()|\rangle, \langle |()\rangle\}.
\end{eqnarray*}
Let us call these strings $b^1,b^2,b^3, b^4$, respectively. Then it follows from the definition that
\begin{eqnarray*}
f_{b^1}(c)&=&\{x_1,x_3x_4-c\},\\
f_{b^2}(c)&=&\{x_1x_2-c,x_4\},\\
f_{b^3}(c)&=&\{x_2, (x_1+x_3)x_4-c\}, \\
f_{b^4}(c)&=&\{x_3,x_1(x_2+x_4)-c\}.
\end{eqnarray*}
Since $u[1,4]=x_1x_2x_3x_4-x_1x_2-x_3x_4-x_1x_4+1$, we see directly that $V(f_{b^j}(c))\subset V_{[1,4]}(1-c)\hspace{1mm}(1\leq j\leq 4)$. Now assume that $n\geq 3$ and that the assertions is proved for any $n'<n$. Take an arbitrary balanced string $b=(b_1,\cdots,b_{2n+1})\in Bra_n$. As is done in the proof of Theorem 4.1, we divide our argument into two cases: (I) when every minimal element in $L(b)$ is maximal, and (II) when there exists an element in $L(b)$ which is minimal but not maximal. \\\\
Case (I) Every minimal element in $L(b)$ is maximal. In this case the left angle bracket must be odd-indexed one $b_{2k-1}$, say, and hence we must have $b_{2k}=|$, and $b_{2k+1}=\rangle$. It follows that
\begin{eqnarray*}
&&cont(b)\\
&&=\{\{1\},\{3\},\cdots,\{2k-3\},\langle 2k-1 | 2k\rangle,\{2k+2\},\{2k+4\},\cdots,\{2n\}\},
\end{eqnarray*}
hence we have
\begin{eqnarray*}
f_b(c)=\{x_1,x_3,\cdots,x_{2k-3},x_{2k-1}x_{2k}-c,x_{2k+2}, x_{2k+4},\cdots,x_{2n}\}.
\end{eqnarray*}
On the other hand the identities (1.6) and (1.7) imply that
\begin{eqnarray*}
&&u(0,x_2,0,x_4,\cdots,0,x_{2k-2},x_{2k-1},x_{2k},x_{2k+1},0, x_{2k+3},\cdots,0, x_{2n-1},0)\\
&=&(-1)^{k-1}\cdot (-1)^{n-k}\cdot u(x_{2k-1},x_{2k})\\
&=&(-1)^{n-1}(x_{2k-1}x_{2k}-1)
\end{eqnarray*}
Therefore if $a=(a_1,\cdots,a_{2n})\in V(f_b(c))$, then we have
\begin{eqnarray*}
u(a_1,\cdots, a_{2n})=(-1)^{n-1}(c-1),
\end{eqnarray*}
and hence we have $a\in V_{[1,2n]}((-1)^{n-1}(c-1))$. \\

\noindent
Case (II): There exists an element in $L(b)$ which is minimal but not maximal. Therefore we are in the situation treated in Proposition 3.2. As is done there, we call such an element $b_k$. Note that this is not a left angle bracket, since the latter is always maximal by the definition of balancedness. Let $b'$ denote the string defined by
\begin{eqnarray*}
b'=(b_1,\cdots,b_{k-1},b_{k+2},\cdots,b_{2n}).
\end{eqnarray*}
Note that  $b'$ belongs to $Bra_{n-1}$. Let $p:\mathbb{A}^{2n}\rightarrow\mathbb{A}^{2n-2}$ denote the map defined by
\begin{eqnarray*}
p(x_1,\cdots,x_{2n})=(x_1,\cdots,x_{k-2},x_{k-1}+x_{k+1},x_{k+2},\cdots,x_{2n}),
\end{eqnarray*}
and let $p_k$ denote its restriction to $V(x_k)\subset\mathbb{A}^{2n}$. Then by a similar reasoning to that employed in deduction of (4.2) and (4.3), we have
\begin{eqnarray}
p_k^{-1}(V(f_{b'}(c)))&=&V(f_b(c)),\\
p_k^{-1}(V_{[1,2n-2]}((-1)^{n-2}(c-1))&\subset & V(x_k)\cap V_{[1,2n]}((-1)^{n-1}(c-1)).
\end{eqnarray}
Since the induction hypothesis implies $V(f_{b'}(c))\subset V_{[1,2n-2]}((-1)^{n-2}(c-1))$, it follows from (5.3) and (5.4) that
\begin{eqnarray*}
V(f_b(c))\subset V(x_k)\cap V_{[1,2n]}((-1)^{n-1}(c-1))\subset V_{[1,2n]}((-1)^{n-1}(c-1)).
\end{eqnarray*}
This completes the proof for Case (II), hence at the same time the proof of Theorem 5.1. \qed\\\\

When $c=1$, the variety $V_{[1,2n]}((-1)^{n-1}(c-1))$ coincides with the original Chebyshev variety $V_{2n}$, we have the following:
\begin{cor}
For any $b\in Bra_n$, we have $V(f_b(1))\subset V_{2n}$.
\end{cor}

\section{Associative transformation}
In this section we investigate when the dimension of the intersection of subvarieties constructed in the previous sections becomes as large as posslble. This will enable us to construct a family of subvarieties of large dimension of AC-varieties. For this purpose we introduce some notation. Firstly we need {\it translation-by-one operator}. For any string $b\in Par_n$ and for any $b_i\in L(b)$, let
\begin{eqnarray*}
cont(b_i)^{(+1)}=\{k+1;k\in cont(b_i)\},
\end{eqnarray*}
and let
\begin{eqnarray*}
cont(b)^{(+1)}=\{cont(b_i)^{(+1)};b_i\in L(b)\}.
\end{eqnarray*}
Accordingly we define the polynomial $f_{b_i}^{(+1)}$ for $b_i\in L(b)$ by
\begin{eqnarray*}
f_{b_i}^{(+1)}=\sum_{p\in cont{b_i}^{(+1)}}x_p,
\end{eqnarray*}
and let
\begin{eqnarray*}
f_b^{(+1)}=\{f_{b_i}^{(+1)};b_i\in L(b)\}.
\end{eqnarray*}
Secondly we need the {\it associative transformation} which is defined as follows:
\begin{df}
For any $b\in Par_n$, suppose that $r(b_1)=b_i$. Let $b'=(b'_j)_{1\leq j\leq 2n}$ be the string of bracket defined by
\begin{displaymath}
b'_j=
\left\{
\begin{array}{ll}
b_{j+1}, & 1\leq j \leq i-2, \\
(, & j =i-1\\ 
b_{j+1}, & i\leq j \leq 2n-1,\\
), & j=2n.
\end{array}
\right.
\end{displaymath}
The string $b'$ is denoted by $ass(b)$, and is called the {\rm associative transformation} of $b$.
\end{df}
\noindent
Remark. The name comes from the usual associative law. If we take the equality $((a+b)+c)+d=(a+b)+(c+d)$ for example, which is legitimate by the associative law and forget the content, then both hand sides become $(())$ and $()()$. The latter is exactly the associative transformation of the former.\\\\
\noindent
Example 6.1. For the string $b=(())()(())\in Par_{5}$ in Example 3.2, the transformed string $b'=ass(b)$ is equal to $b'=()(()(()))$. Therefore the set of polynomials $f_{ass(b)}^{(+1)}$ is given by
\begin{eqnarray*}
f_{ass(b)}^{(+1)}=\{x_2,x_5,x_7+x_9,x_8,x_4+x_6+x_{10}\}
\end{eqnarray*}
Note here there is a strong resemblance between this and the set 
\begin{eqnarray*}
f_b=\{x_1+x_3,x_2,x_5,x_7+x_9,x_8\}.
\end{eqnarray*}
The elements of $f_{ass(b)}^{(+1)}$ except the last one and the elements of $f_b$ except the first one coincide completely. This phenomina is a special case of the following general proposition;

\begin{prp}
For any $b\in Par_n$, its associative transformation $ass(b)$ is balanced too, hence it defines a map $ass:Par_n\rightarrow Par_n$. Moreover we have
\begin{eqnarray*}
f_b\cap K[x_2,\cdots,x_{2n-1}]=f_{ass(b)}^{(+1)}\cap K[x_2,\cdots,x_{2n-1}],
\end{eqnarray*}
and
\begin{eqnarray*}
&&\#(f_b\setminus (f_b\cap K[x_2,\cdots,x_{2n-1}]))\\
&&=\#(f_{ass(b)}^{(+1)}\setminus (f_{ass(b)}^{(+1)}\cap K[x_2,\cdots,x_{2n-1}]))=1.
\end{eqnarray*}
Therefore we have
\begin{eqnarray*}
\#(f_b\cup f_{ass(b)}^{(+1)})=n+1.
\end{eqnarray*}
\end{prp}

\noindent
{\it Proof}. This is a direct consequence of the definition of the content and the associative transformation. \qed\\\\
This has a nice consequence for the geometry of the intersection of two Chebyshev varieties:

\begin{thm}
For any integer $n\geq2$ and for any $b\in Par_n$, we have
\begin{eqnarray}
V(f_b)\cap V(f_{ass(b)}^{(+1)})\subset V(u[1,2n-1])\cap V(u[2,2n]),
\end{eqnarray}
and the dimension of each side is given by
\begin{eqnarray}
\dim (V(f_b)\cap V(f_{ass(b)}^{(+1)}))&=&n-1,\\
\dim (V(u[1,2n-1])\cap V(u[2,2n]))&=&2(n-1).
\end{eqnarray}
\end{thm}

\noindent
{\it Proof}. The assertion (6.1) follows from Corrolary 4.1. As for (6.2), recall that $cont(b_i)\cap cont(b_j)=\phi$ for any pair of distinct left brackets $b_i, b_j\in L(b)$ by the very definition of the content. Therefore the linear polynomials $f_{b_i}\hspace{1mm}(b_i\in L(b))$ are linearly independent. Moreover an element of $f_b$ contains the term $x_1$, and an element of $f_{ass(b)}^{(+1)}$ contains the term $x_{2n}$, as is noted in Proposition 6.1. Hence the $n+1$ elements of $f_b\cup f_{ass(b)}^{(+1)}$ are linearly independent and the assertion (6.2) follows. The last assertion (6.3) is a direct consequence of the fact that AC-variety is a complete intersection. This completes the proof. \qed\\\\

\section{Bra-ketting transformation}
In order to construct subvarieties of AC-varieties and to require the dimension of those as large as possible, we introduce another type of transformation of the string of round brackets, called {\it bra-ketiing transformation}.

\begin{df}
For any $b\in Par_n$ with $r(b_1)=b_i$, let $b'=(b'_j)_{1\leq j\leq 2n+1}$ be the string of bracket defined by
\begin{displaymath}
b'_j=
\left\{
\begin{array}{ll}
\langle, & j=1, \\
b_j, & 2\leq j\leq i-1,\\ 
|, & j=i,\\
b_j, & i+1\leq j\leq 2n,\\
\rangle, & j=2n+1.
\end{array}
\right.
\end{displaymath}
The string $b'$ is denoted by $bra(b)$, and is called the {\rm bra-ketting transformation} of $b$.
\end{df}

\noindent
Example 7.1. For the string $b=(())()(())\in Par_{5}$ in Example 3.2, the transformed string $b'=bra(b)$ is equal to $b'=\langle()|()(())\rangle$. By the definition, for any $b\in Par_n$, the bra-ketting transformation $bra(b)$ is balanced. Hence it defines a map $bra:Par_n\rightarrow Bra_n$. This transformation enables us to construct a linear subvarieties of $AC_N$ when $N\equiv 0\pmod 4$:

\begin{thm}
For any positive integer $n$ and for any $b\in Par_n$, we have
\begin{eqnarray}
V(f_b)\cap V(f_{ass(b)}^{(+1)})\subset V(f_{bra(b)}(0)).
\end{eqnarray}
\end{thm}

\noindent
{\it Proof}. Let $b'=ass(b)$, $d=bra(b)$. Let $L(b)=\{b_{\ell_j};1\leq j\leq n\}$ with $\ell_1<\cdots <\ell_n$ so that $\ell_1=1$, and let $r(b_1)=b_i$. It follows from the definition that
\begin{eqnarray*}
cont(d_1)&=&cont(b_1),\\
cont(d_i)&=&cont(b'_{i-1})^{(+1)}.
\end{eqnarray*}
Moreover, for any left bracket $b_{\ell_j}$ with $j\geq 2$, we have
\begin{eqnarray*}
cont(b_{\ell_j})=cont(b'_{\ell_j-1})^{(+1)}=cont(d_{\ell_j}).
\end{eqnarray*}
(Note here that $b_i\notin L(b)$.) Therefore we have
\begin{eqnarray*}
f_b\cup f_{b'}^{(+1)}=\{f_{b_1}\}\cup \{f_{b'_{i-1}}^{(+1)}\}\cup \{f_{b_{\ell_j}};2\leq j\leq n\},
\end{eqnarray*}
and
\begin{eqnarray*}
f_d(0)&=&\{f_{b_{\ell_j}};2\leq j\leq n\}\cup \{f_{d_1}(0)\}\\
&=&\{f_{b_{\ell_j}};2\leq j\leq n\}\cup \{f_{b_1}\cdot f_{b'_{i-1}}^{(+1)}\}.
\end{eqnarray*}
Hence the inclusion (7.1) holds true. This completes the proof. \qed\\\\

Recall that for any $d\in Bra_n$ and for any $c\in K$, we have $V(f_d(c))\subset V_{[1,2n]}((-1)^{n-1}(c-1))$ by Theorem 6.1. Hence when $n$ is even, we have
\begin{eqnarray}
V(f_d(0))\subset V_{[1,2n]}((-1)^n)=V_{[1,2n]}(1).
\end{eqnarray}
Now the AC-variety $AC_N$ is defined to be the intersection of three hypersurfaces $V_{[1,N-1]}$, $V_{[2,N]}$, and $V_{[1,N]}(1)$ in $\mathbb{A}^N$. On the other hand, for any $b\in Par_{2n}$, it follows from Corollary 4.1 that
\begin{eqnarray}
V(f_b)&\subset& V_{[1,2n-1]},\\
V(f_{ass(b)}^{(+1)})&\subset& V_{[2,2n]}.
\end{eqnarray}
Hence combining the inclusions (7.1)-(7.4), we arrive at the following:

\begin{thm}
For any positive integer $n\equiv 0 \pmod 4$, and for any $b\in Par_{n/2}$, we have
\begin{eqnarray*}
V(f_b)\cap V(f_{ass(b)}^{(+1)})\subset AC_n.
\end{eqnarray*}
\end{thm}

\section{Strings with a triple bra-ket}
Here a {\it triple bra-ket} means the string $\langle\hspace{1mm}|\hspace{1mm}|\hspace{1mm}\rangle$. This will enable us to construct subvarieties of $AC_N$ for an odd integer $N$. For any integer $n\geq 1$, we define the balancedness of a string $b=(b_1,\cdots, b_{2n+2})$ with $n-1$ pairs of round brackets and with one triple bra-ket as follows. Let $b_i=\langle, b_j=|, b_k=|, b_{\ell}=\rangle$ with $i<j<k<\ell$. Then $b$ is said to be balanced if the five remaining substrings$(b_1,\cdots,b_{i-1}), (b_{i+1},\cdots, b_{j-1}), (b_{j+1},\cdots, b_{k-1})$, $(b_{k+1},\cdots, b_{\ell-1})$, and $(b_{\ell+1},\cdots,b_{2n+2})$ are balanced in the sense of Definition 3.1. Moreover we employ the convention that the left angle bracket $b_i$ matches the first vertical bar $b_j$, which matches the second vertical bar $b_k$, which matches the right angle bracket $b_{\ell}$. Accordingly we regard a vertical bar both as a left bracket and as a right bracket, and we can define a partial order on the set $L(b)$ of left brackets in $b$ in the same way as for $Par_n$. We denote the set of balanced strings of length $2n+2$ with $n-1$ pairs of round brackets and with one triple bra-ket by $Tbra_n$. The content of a string $b\in Tbra_n$ is defined in the same way as for $Par_n$ with the only difference that the content of the triple bra-ket $(b_i,b_j,b_k,b_{\ell})$ is denoted by $\langle cont(b_i) | cont(b_j)| cont(b_k) \rangle$. For $\epsilon\in \{\pm 1\}$, we associate to $b_i, b_j, b_k$ the polynomials $f_{b_i}(\epsilon), f_{b_j}(\epsilon), f_{b_k}(\epsilon)$ defined by
\begin{eqnarray*}
f_{b_i}(\epsilon)&=&\sum_{p\in cont(b_i)}x_p-\epsilon,\\
f_{b_j}(\epsilon)&=&\sum_{p\in cont(b_j)}x_p-\epsilon,\\
f_{b_k}(\epsilon)&=&\sum_{p\in cont(b_k)}x_p-\epsilon.
\end{eqnarray*}
For the remaining left brackets $b_m, m\neq i,j,k$, the polynomial $f_{b_m}$ is defined to be the same as in the case of round bracket. Putting them together, we define the set of polynomials $f_b(\epsilon)$ by
\begin{eqnarray*}
f_b(\epsilon)=\{f_{b_m};1\leq m\leq n, m\neq i,j,k\}\cup \{f_{b_i}(\epsilon), f_{b_j}(\epsilon), f_{b_k}(\epsilon)\}.
\end{eqnarray*}
The following example should clarify the meaning.\\\\

\noindent
Example 8.1. Let $b=()\langle(())|()|\rangle\in Tbra_{5}$. Then the content and the associated polynomial are given by
\begin{eqnarray*}
cont(b)&=&\{{1},{5},{9},\{4,6\},\langle 3,7 | 8,10 | 11\rangle\},\\
f_b(\epsilon)&=&\{x_1, x_5, x_9, x_4+x_6, x_3+x_7-\epsilon, x_8+x_{10}-\epsilon, x_{11}-\epsilon\}.
\end{eqnarray*}
\\
The family of polynomials $f_b(\epsilon)$ for $b\in Tbra_n$ will play an important role to construct subvarieties of AC-varieties. First we prove the following:
\begin{prp}
For any integer $n\geq 1$, $c\in K$, and for any $b\in Tbra_n$, we have
\begin{eqnarray*}
V(f_b(\epsilon))\subset V_{[1,2n+1]}((-1)^n\epsilon).
\end{eqnarray*}
\end{prp}

\noindent
{\it Proof}.
We prove this by induction on $n$. When $n=1$, the set $Tbra_1$ consists solely of $\langle\hspace{1mm}|\hspace{1mm}|\hspace{1mm}\rangle$. If we call this string $b$, then its content is given by
\begin{eqnarray*}
cont(b)=\{\langle1|2|3\rangle\}.
\end{eqnarray*}
Therefore we have
\begin{eqnarray*}
f_b(\epsilon)=\{x_1-\epsilon, x_2-\epsilon, x_3-\epsilon\}.
\end{eqnarray*}
On the other hand, we have 
\begin{eqnarray*}
V_{[1,3]}((-1)^1\epsilon)&=&\{(x_1,x_2,x_3)\in\mathbb{A}^3;u[1,3]=-\epsilon\}\\
&=&\{(x_1,x_2,x_3)\in\mathbb{A}^3;x_1x_2x_3-x_1-x_3=-\epsilon\},
\end{eqnarray*}
hence the assertion holds since $\epsilon^3=\epsilon$ for $\epsilon\in \{\pm 1\}$. When $n=2$, the set $Tbra_2$ has five strings:
\begin{eqnarray*}
Tbra_2=\{()\langle ||\rangle, \langle ()||\rangle, \langle |()|\rangle, \langle ||()\rangle, \langle ||\rangle ()\}.
\end{eqnarray*}
Let us call these strings $b^1,b^2,b^3, b^4, b^5$, respectively. Then it follows from the definition that
\begin{eqnarray*}
f_{b^1}(\epsilon)&=&\{x_1,x_3-\epsilon,x_4-\epsilon, x_5-\epsilon\},\\
f_{b^2}(\epsilon)&=&\{x_2,x_1+x_3-\epsilon, x_4-\epsilon,x_5-\epsilon\},\\
f_{b^3}(\epsilon)&=&\{x_3, x_1-\epsilon, x_2+x_4-\epsilon,x_5-\epsilon\}, \\
f_{b^4}(\epsilon)&=&\{x_4,x_1-\epsilon, x_2-\epsilon, x_3+x_5-\epsilon\},\\
f_{b^5}(\epsilon)&=&\{x_5,x_1-\epsilon, x_2-\epsilon, x_3-\epsilon\}.
\end{eqnarray*}
Since $u[1,5]=x_1x_2x_3x_4x_5-x_1x_2x_3-x_1x_2x_5-x_1x_4x_5-x_3x_4x_5+x_1+x_3+x_5$, we see directly that $V(f_{b^j}(\epsilon))\subset V_{[1,5]}(\epsilon)\hspace{1mm}(1\leq j\leq 5)$. Now assume that $n\geq 3$ and that the assertions is proved for any $n'<n$. Take an arbitrary balanced string $b=(b_1,\cdots,b_{2n+2})\in Tbra_n$. As is done in the proof of Theorem 4.1, we divide our argument into two cases: (I) when every minimal element in $L(b)$ is maximal, and (II) when there exists an element in $L(b)$ which is minimal but not maximal. \\\\
Case (I) Every minimal element in $L(b)$ is maximal. In this case the left angle bracket must be odd-indexed one $b_{2k-1}$, say, and hence we must have $b_{2k}=|$, $b_{2k+1}=|$, and $b_{2k+2}=\rangle$. It follows that
\begin{eqnarray*}
&&cont(b)\\
&&=\{\{1\},\{3\},\cdots,\{2k-3\},\langle 2k-1 | 2k | 2k+1 \rangle,\{2k+3\},\cdots,\{2n+1\}\},
\end{eqnarray*}
hence we have
\begin{eqnarray*}
&&f_b(\epsilon)\\
&&=\{x_1,x_3,\cdots,x_{2k-3},x_{2k-1}-\epsilon,x_{2k}-\epsilon,x_{2k+1}-\epsilon,x_{2k+3}, \cdots,x_{2n+1}\}.
\end{eqnarray*}
On the other hand the equalities (1.6) and (1.7) imply that
\begin{eqnarray*}
&&u(0,x_2,0,\cdots,0,x_{2k-2},x_{2k-1},x_{2k},x_{2k+1},x_{2k+2},0, x_{2k+4},\cdots,0, x_{2n},0)\\
&=&(-1)^{k-1}\cdot (-1)^{n-k}\cdot u(x_{2k-1},x_{2k},x_{2k+1})\\
&=&(-1)^{n-1}(x_{2k-1}x_{2k}x_{2k+1}-x_{2k-1}-x_{2k+1}).
\end{eqnarray*}
Therefore if $a=(a_1,\cdots,a_{2n+1})\in V(f_b(\epsilon))$, then we have
\begin{eqnarray*}
u(a_1,\cdots, a_{2n+1})=(-1)^{n-1}\cdot(-\epsilon)=(-1)^n\epsilon,
\end{eqnarray*}
and hence we have $a\in V_{[1,2n+1]}((-1)^n\epsilon)$. \\

\noindent
Case (II): There exists an element in $L(b)$ which is minimal but not maximal. We call such an element $b_k$. Note that this is not a left angle bracket nor a vertical bar, since these are always maximal by the definition of balancedness. Let $b'$ denote the string defined by
\begin{eqnarray*}
b'=(b_1,\cdots,b_{k-1},b_{k+2},\cdots,b_{2n+1}).
\end{eqnarray*}
Note that  $b'$ belongs to $Tbra_{n-1}$. Let $p:\mathbb{A}^{2n+1}\rightarrow\mathbb{A}^{2n-	1}$ denote the map defined by
\begin{eqnarray*}
p(x_1,\cdots,x_{2n+1})=(x_1,\cdots,x_{k-2},x_{k-1}+x_{k+1},x_{k+2},\cdots,x_{2n+1}),
\end{eqnarray*}
and let $p_k$ denote its restriction to $V(x_k)\subset\mathbb{A}^{2n+1}$. Then by a similar reasoning to that for (4.2) and (4.2) we see that
\begin{eqnarray}
p_k^{-1}(V(f_{b'}(\epsilon)))&=&V(f_b(\epsilon)),\\
p_k^{-1}(V_{[1,2n-1]}((-1)^{n-1}\epsilon)&\subset & V(x_k)\cap V_{[1,2n+1]}((-1)^n\epsilon).
\end{eqnarray}
Since the induction hypothesis implies $V(f_{b'}(\epsilon))\subset V_{[1,2n-1]}((-1)^{n-1}\epsilon)$, it follows from (8.1) and (8.2) that
\begin{eqnarray*}
V(f_b(\epsilon))\subset V(x_k)\cap V_{[1,2n+1]}((-1)^n\epsilon)\subset V_{[1,2n+1]}((-1)^n\epsilon).
\end{eqnarray*}
This completes the proof for Case (II), hence at the same time the proof of Proposition 8.1. \qed\\\\

\section{Associative transformation II}
In this section we introduce another kind of associative transformation. This will enable us to construct a family of subvarieties of $AC_N$ for odd $N$. Let $Bra_{\langle n}$ denote the subset of $Bra_n$ consisting of strings with $b_1=\langle$. We define {\it associative transformation for bra-ket} as follows:
\begin{df}
For any $b\in Bra_{\langle n}$, suppose that $b_i=|$ and $b_j=\rangle$. Let $b'=(b'_k)_{1\leq k\leq 2n+1}$ be the string of bracket defined by
\begin{displaymath}
b'_k=
\left\{
\begin{array}{ll}
b_{k+1}, & 1\leq k \leq i-2, \\
\langle, & k =i-1\\ 
b_{k+1}, & i\leq k \leq j-2,\\
|, & k=j-1,\\
b_{k+1}, & j\leq k\leq 2n,\\
\rangle, & k=2n+1.
\end{array}
\right.
\end{displaymath}
The string $b'$ is denoted by $ass_{\langle |\rangle}(b)$, and is called the {\rm associative transformation for bra-ket} of $b$.
\end{df}
\noindent
Example 9.1. Let $b=\langle()|(())\rangle ()\in Bra_{5}$. Then the transformed string $b'=ass_{\langle|\rangle}(b)$ is equal to $b'=()\langle(())|()\rangle$. Therefore $f_b(c)$ and $f_{ass_{\langle|\rangle}(b)}^{(+1)}(c)$ are given by
\begin{eqnarray*}
f_b(c)=\{x_2,x_6,x_5+x_7,x_{10},(x_1+x_3)(x_4+x_8)-c\},
\end{eqnarray*}
\begin{eqnarray*}
f_{ass_{\langle|\rangle}(b)}^{(+1)}(c)=\{x_2,x_6,x_5+x_7,x_{10},(x_4+x_8)(x_9+x_{11})-c\}.
\end{eqnarray*}
Note here too that there is a strong resemblance between $f_b(c)$ and $f_{ass_{\langle|\rangle}(b)}^{(+1)}(c)$, and we have the following proposition which can be proved in the same way as for Proposition 6.1:

\begin{prp}
For any $b\in Bra_n$, its associative transformation for bra-ket $ass_{\langle|\rangle}(b)$ is balanced too, hence it defines a map $ass_{\langle|\rangle}:Bra_n\rightarrow Bra_n$. Moreover we have
\begin{eqnarray*}
f_b(c)\cap K[x_2,\cdots,x_{2n}]=f_{ass(b)}^{(+1)}(c)\cap K[x_2,\cdots,x_{2n}],
\end{eqnarray*}
and
\begin{eqnarray*}
&&\#(f_b(c)\setminus (f_b(c)\cap K[x_2,\cdots,x_{2n-1}]))\\
&&=\#(f_{ass(b)}^{(+1)}(c)\setminus (f_{ass(b)}^{(+1)}(c)\cap K[x_2,\cdots,x_{2n-1}]))\\
&&=1.
\end{eqnarray*}
Therefore we have
\begin{eqnarray*}
\#(f_b(c)\cup f_{ass(b)}^{(+1)}(c))=n-1.
\end{eqnarray*}
\end{prp}

\noindent
When we set $c$ to be $\epsilon\in \{\pm 1\}$ in the pair $f_b(c)$ and $f_{ass_{\langle|\rangle}(b)}^{(+1)}(c)$, we find a strong connection between these and the set of polynomial associated to a string $\in Tbra_n$. In order to explain the connection we introduce one more transformation:

\begin{df}
For any $b\in Bra_{\langle n}$, suppose that $b_i=|$ and $b_j=\rangle$. Let $b"=(b"_k)_{1\leq k\leq 2n+2}$ be the string of bracket defined by
\begin{displaymath}
b"_k=
\left\{
\begin{array}{ll}
\langle, & k=1,\\
b_k, & 2\leq k \leq i-1, \\
|, & k =i\\ 
b_k, & i+1\leq k \leq j-1,\\
|, & k=j,\\
b_k, & j+1\leq k\leq 2n+1,\\
\rangle, & k=2n+2.
\end{array}
\right.
\end{displaymath}
The string $b"$ is denoted by $ass_{to\langle ||\rangle}(b)$, and is called the {\rm associative transformation to triple bra-ket} of $b$.
\end{df}

\noindent
Example 9.2. For the string $b=\langle()|(())\rangle ()\in Bra_{5}$ used in Example 9.1, the transformed string $b"=ass_{to\langle||\rangle}(b)$ is equal to $\langle()|(())|()\rangle$. The following proposition is crucial for our construction of a family of linear subvarieties of $AC_N$ with $N$ odd:

\begin{prp}
For any integer $n\geq 1$ and for any string $b\in Bra_{\langle n}$, let $b'=ass_{\langle |\rangle}(b)$ and $b"=ass_{to\langle ||\rangle}(b)$. Then for any $\epsilon\in \{\pm 1\}$, we have
\begin{eqnarray}
V(f_{b"}(\epsilon))\subset V(f_b(1))\cap V(f_{b'}^{(+1)}(1)).
\end{eqnarray}
\end{prp}

\noindent
{\it Proof}. Let $b_i=|$ and $b_j=\rangle$. Note that $b_1=\langle$ by the definition of $Bra_{\langle n}$. Then it follows from the definition that for any $b_k\in L(b)-\{b_1,b_i\}$, we have 
\begin{eqnarray*}
f_{b_k}=f_{b'_{k-1}}^{(+1)}=f_{b"_k}.
\end{eqnarray*}
Therefore we have
\begin{eqnarray}
f_b(1)&=&\{f_{b_k};b_k\in L(b)-\{b_1,b_i\}\}\cup \{f_{b_1}(1)\},\\
f_{b'}^{(+1)}(1)&=&\{f_{b_k};b_k\in L(b)-\{b_1,b_i\}\}\cup \{f_{b'_{i-1}}^{(+1)}(1)\},\\
f_{b"}(\epsilon)&=&\{f_{b_k};b_k\in L(b)-\{b_1,b_i\}\}\cup \{f_{b"_1}(\epsilon), f_{b"_i}(\epsilon), f_{b"_j}(\epsilon) \}.
\end{eqnarray}
Recall that the polynomials $f_{b_1}(1),\cdots,  f_{b"_j}(\epsilon)$ on the right hand sides of the above equalities are defined as
\begin{eqnarray*}
f_{b_1}(1)&=&\sum_{p\in cont(b_1)}x_p\sum_{q\in cont(b_i)}x_q-1,\\
f_{b'_{i-1}}^{(+1)}(1)&=&\sum_{q\in cont(b'_{i-1})}x_{q+1}\sum_{r\in cont(b'_{j-1})}x_{r+1}-1\\
&=&\sum_{q\in cont(b_i)}x_q\sum_{r\in cont(b'_{j-1})}x_{r+1}-1,\\
f_{b"_1}(\epsilon)&=&\sum_{p\in cont(b"_1)}x_p-\epsilon\\
&=&\sum_{p\in cont(b_1)}x_p-\epsilon,\\
f_{b"_i}(\epsilon)&=&\sum_{q\in cont(b"_i)}x_q-\epsilon\\
&=&\sum_{q\in cont(b_i)}x_q-\epsilon,\\
f_{b"_j}(\epsilon)&=&\sum_{r\in cont(b"_j)}x_r-\epsilon\\
&=&\sum_{r\in cont(b'_{j-1})}x_{r+1}-\epsilon.
\end{eqnarray*}
Combining these with (9.2)-(9.4) and noting that $\epsilon^2=1$ for  $\epsilon\in \{\pm 1\}$, we see that the inclusion (9.1) holds true. This completes the proof. \qed\\\\

\noindent
This proposition enables us to obtain a family of linear subvarieties of $AC_N$ for an arbitrary odd integer $N$:

\begin{thm}
For any integer $n\geq 1$ and for any string $b\in Bra_{\langle n}$, we have
\begin{eqnarray}
V(f_{ass_{to\langle||\rangle}(b)}((-1)^n))\subset AC_{2n+1}.
\end{eqnarray}
\end{thm}

\noindent
{\it Proof}. It follows from Corollary 5.1 that
\begin{eqnarray}
V(f_b(1))&\subset& V_{[1,2n]},\\
V(f_{ass_{\langle|\rangle}(b)}^{(+1)}(1))&\subset& V_{[2,2n+1]}.
\end{eqnarray}
Furthermore it follows from Proposition 8.1 that
\begin{eqnarray}
V(f_{ass_{to\langle||\rangle}(b)}((-1)^n))\subset V_{[1,2n+1]}((-1)^n\cdot (-1)^n)=V_{[1,2n+1]}(1).
\end{eqnarray}
Hence combining (9.6)-(9.8) with Proposition 9.2, we have the inclusion (9.5). This completes the proof. \qed\\\\

\section{Strings with a quadruple bra-ket}
In order to deal with subvarieties of $AC_N$ with $N\equiv 2\pmod 4$, we need strings with a {\it quadruple bra-ket}. Here a quadruple bra-ket means the string $\langle\hspace{1mm}|\hspace{1mm}|\hspace{1mm}|\hspace{1mm}\rangle$. For any integer $n\geq 2$, we define the balancedness of a string $b=(b_1,\cdots, b_{2n+1})$ with $n-2$ pairs of round brackets and with one quadruple bra-ket in a way similar to that for strings with a triple bra-ket. Let $b_i=\langle, b_j=|, b_k=|, b_{\ell}=|, b_m=\rangle$ with $i<j<k<\ell<m$. Then $b$ is said to be balanced if the six remaining substrings$(b_1,\cdots,b_{i-1}), (b_{i+1},\cdots, b_{j-1}), (b_{j+1},\cdots, b_{k-1})$, $(b_{k+1},\cdots, b_{\ell-1})$, $(b_{\ell +1},\cdots, b_{m-1})$, and $(b_{m+1},\cdots,b_{2n+1})$ are balanced in the sense of (). We regard a vertical bar both as a left bracket and as a right bracket, and we can define a partial order on the set $L(b)$ of left brackets in $b$ in the same way as for $Par_n$. We denote the set of balanced strings of length $2n+1$ with $n-2$ pairs of round brackets and with one quadruple bra-ket by $Qbra_n$. The content of a string $b\in Qbra_n$ is defined in the same way as for $Par_n$ with the only difference that the content of the bra-ket $(b_i,b_j,b_k,b_{\ell}, b_m)$ is denoted by $\langle cont(b_i) | cont(b_j)| cont(b_k)| cont(b_{\ell}) \rangle$. We associate to the vertical bars $b_j, b_k, b_{\ell},$ the polynomials $f_{b_j}(2), f_{b_k}(2), f_{b_{\ell}}(2)$ defined by
\begin{eqnarray*}
f_{b_j}(2)&=&\sum_{p\in cont(b_i)}x_p\sum_{q\in cont(b_j)}x_q-2,\\
f_{b_k}(2)&=&\sum_{p\in cont(b_j)}x_p\sum_{q\in cont(b_k)}x_q-2,\\
f_{b_{\ell}}(2)&=&\sum_{p\in cont(b_k)}x_p\sum_{q\in cont(b_{\ell})}x_q-2.
\end{eqnarray*}
For the remaining left brackets $b_r, r\neq i,j,k,\ell$, the polynomial $f_{b_r}$ is defined to be the same as in the case of round bracket. Putting them together, we define the set of polynomials $f_b(2)$ by
\begin{eqnarray*}
f_b(2)=\{f_{b_m};b_m\in L(b), m\neq i,j,k,\ell\}\cup \{f_{b_j}(2), f_{b_k}(2), f_{b_l}(2)\}.
\end{eqnarray*}
\noindent
The family of polynomials $f_b(2)$ for $b\in Qbra_n$ will play an important role to construct subvarieties of AC-varieties $AC_N$ with $N\equiv 2\pmod 4$. First we prove the following:

\begin{prp}
For any integer $n\geq 2$ and for any $b\in Qbra_n$, we have
\begin{eqnarray*}
V(f_b(2))\subset V_{[1,2n]}((-1)^{n-1}).
\end{eqnarray*}
\end{prp}

\noindent
{\it Proof}.
We prove this by induction on $n$. When $n=2$, the set $Qbra_2$ consists solely of $\langle\hspace{1mm}|\hspace{1mm}|\hspace{1mm}|\hspace{1mm}\rangle$. If we call this string $b$, then its content is given by
\begin{eqnarray*}
cont(b)=\{\langle1|2|3|4\rangle\}.
\end{eqnarray*}
Therefore we have
\begin{eqnarray*}
f_b(2)=\{x_1x_2-2, x_2x_3-2, x_3x_4-2\}.
\end{eqnarray*}
Note that if $(a_1,\cdots, a_4)\in V(f_b(2))$, then $a_i\neq 0$ for $i=1, \cdots, 4$, hence the second and the third equations imply that $a_2=a_4$ holds. On the other hand, we have 
\begin{eqnarray*}
&&V_{[1,4]}((-1)^{2-1})\\
&=&\{(x_1,x_2,x_3, x_4)\in\mathbb{A}^4;u[1,4]=-1\}\\
&=&\{(x_1,x_2,x_3,x_4)\in\mathbb{A}^4;x_1x_2x_3x_4-x_1x_2-x_1x_4-x_3x_4+1=-1\},
\end{eqnarray*}
hence if $(a_1,\cdots, a_4)\in V(f_b(2))$, then we have
\begin{eqnarray*}
a_1a_2a_3a_4-a_1a_2-a_1a_4-a_3a_4+1=2\cdot 2-2-2-2+1=-1,
\end{eqnarray*}
from which the assertion follows. Now assume that $n\geq 3$ and that the assertions is proved for any $n'<n$. Take an arbitrary balanced string $b=(b_1,\cdots,b_{2n+1})\in Qbra_n$. As is done in the proof of Theorem 4.1, we divide our argument into two cases: (I) when every minimal element in $L(b)$ is maximal, and (II) when there exists an element in $L(b)$ which is minimal but not maximal. \\\\
Case (I) Every minimal element in $L(b)$ is maximal. In this case the left angle bracket must be odd-indexed one $b_{2k-1}$, say, and hence we must have $b_{2k}=|$, $b_{2k+1}=|$, $b_{2k+2}=|$, and $b_{2k+3}=\rangle$. It follows that
\begin{eqnarray*}
&&cont(b)\\
&=&\{\{1\},\{3\},\cdots,\{2k-3\},\langle 2k-1 | 2k | 2k+1 | 2k+2 \rangle,\{2k+4\},\cdots,\{2n\}\},
\end{eqnarray*}
hence we have
\begin{eqnarray*}
f_b(2)&=&\{x_1,x_3,\cdots,x_{2k-3},x_{2k-1}x_{2k}-2,x_{2k}x_{2k+1}-2,\\
&&\hspace{8em}x_{2k+1}x_{2k+2}-2,x_{2k+4}, \cdots,x_{2n}\}.
\end{eqnarray*}
On the other hand the equalities (1.6) and (1.7) imply that
\begin{eqnarray*}
&&u(0,x_2,0,\cdots,0,x_{2k-2},x_{2k-1},x_{2k},x_{2k+1},x_{2k+2},x_{2k+3},\\
&&\hspace{16em}0, x_{2k+5},\cdots,0, x_{2n-1},0)\\
&=&(-1)^{k-1}\cdot (-1)^{n-k-1}\cdot u(x_{2k-1},x_{2k},x_{2k+1},x_{2k+2})\\
&=&(-1)^nu(x_{2k-1},x_{2k},x_{2k+1},x_{2k+2}).
\end{eqnarray*}
Therefore if $a=(a_1,\cdots,a_{2n})\in V(f_b(2))$, then it follows from our proof for the case $n=2$ that
\begin{eqnarray*}
u(a_1,\cdots, a_{2n+1})=(-1)^n\cdot(-1)=(-1)^{n-1},
\end{eqnarray*}
and hence we have $a\in V_{[1,2n]}((-1)^{n-1})$. \\

\noindent
Case (II): There exists an element in $L(b)$ which is minimal but not maximal. We call such an element $b_k$. Note that this is not a left angle bracket nor a vertical bar, since those are always maximal by the definition of balancedness. Let $b'$ denote the string defined by
\begin{eqnarray*}
b'=(b_1,\cdots,b_{k-1},b_{k+2},\cdots,b_{2n+1}).
\end{eqnarray*}
Note that  $b'$ belongs to $Qbra_{n-1}$. Let $p:\mathbb{A}^{2n}\rightarrow\mathbb{A}^{2n-2}$ denote the map defined by
\begin{eqnarray*}
p(x_1,\cdots,x_{2n})=(x_1,\cdots,x_{k-2},x_{k-1}+x_{k+1},x_{k+2},\cdots,x_{2n}),
\end{eqnarray*}
and let $p_k$ denote its restriction to $V(x_k)\subset\mathbb{A}^{2n}$. Then by a similar reasoning to that for (4.2) and (4.3), we have
\begin{eqnarray}
p_k^{-1}(V(f_{b'}(2)))&=&V(f_b(2)),\\
p_k^{-1}(V_{[1,2n-2]}((-1)^{n-2})&\subset & V(x_k)\cap V_{[1,2n]}((-1)^{n-1}).
\end{eqnarray}
Since the induction hypothesis implies $V(f_{b'}(2))\subset V_{[1,2n-2]}((-1)^{n-2})$, it follows from (10.1) and (10.2) that
\begin{eqnarray*}
V(f_b(2))\subset V(x_k)\cap V_{[1,2n]}((-1)^{n-1})\subset V_{[1,2n]}((-1)^{n-2}).
\end{eqnarray*}
This completes the proof for Case (II), hence at the same time the proof of Proposition 10.1. \qed\\\\

In order to construct a family of subvarieties of $AC_N$ with $N\equiv 2\pmod 4$, we need to introduce sets of polynomials $f_b(2)$ for $b\in Tbra_n$ too. These are, however, defined in a similar way to those for $b\in Qbra_n$ as follows. For any $b\in Tbra_n$, assume that $b_i=\langle, b_j=|, b_k=|, b_{\ell}=\rangle$ with $i<j<k<\ell$. We associate to the vertical bars $b_j, b_k$ the polynomials $f_{b_j}(2), f_{b_k}(2)$ defined by
\begin{eqnarray*}
f_{b_j}(2)&=&\sum_{p\in cont(b_i)}x_p\sum_{q\in cont(b_j)}x_q-2,\\
f_{b_k}(2)&=&\sum_{p\in cont(b_j)}x_p\sum_{q\in cont(b_k)}x_q-2.\\
\end{eqnarray*}
For the remaining left brackets $b_m, m\neq i,j,k$, the polynomial $f_{b_m}$ is defined to be the same as in the case of round bracket. Putting them together, we define the set of polynomials $f_b(2)$ by
\begin{eqnarray*}
f_b(2)=\{f_{b_m};b_m\in L(b), m\neq i,j,k\}\cup \{f_{b_j}(2), f_{b_k}(2)\}.
\end{eqnarray*}

\noindent
Then we can show the following:

\begin{prp}
For any integer $n\geq 1$ and for any $b\in Tbra_n$, we have
\begin{eqnarray*}
V(f_b(2))\subset V_{[1,2n+1]}.
\end{eqnarray*}
\end{prp}

\noindent
{\it Proof}. We prove this by induction on $n$. Since our proof goes similarly to that for Proposition 10.1, we only point out some necessary minor changes to deal with this case. When $n=1$, the set $Tbra_1$ consists solely of $\langle\hspace{1mm}|\hspace{1mm}|\hspace{1mm}\rangle$. If we call this string $b$, then its content is given by
\begin{eqnarray*}
cont(b)=\{\langle1|2|3\rangle\}.
\end{eqnarray*}
Therefore we have
\begin{eqnarray*}
f_b(2)=\{x_1x_2-2, x_2x_3-2\}.
\end{eqnarray*}
Note that if $(a_1,\cdots, a_3)\in V(f_b(2))$, then $a_i\neq 0$ for $i=1, 2, 3$, hence the first and the second equations imply that $a_1=a_3$. On the other hand, we have 
\begin{eqnarray*}
V_{[1,3]}&=&\{(x_1,x_2,x_3)\in\mathbb{A}^3;u[1,3]=0\}\\
&=&\{(x_1,x_2,x_3)\in\mathbb{A}^3;x_1x_2x_3-x_1-x_3=0\},
\end{eqnarray*}
hence if $(a_1,\cdots, a_3)\in V(f_b(2))$, then we have
\begin{eqnarray*}
a_1a_2a_3-a_1-a_3=2a_3-a_1-a_3=0,
\end{eqnarray*}
from which the assertion holds. The remaining induction process can be proved similarly and we omit it. \qed\\

\section{Associative transformation III}
In this section we introduce one more associative transformation. This will enable us to construct a family of subvarieties of $AC_N$ for $N\equiv 2 \pmod 4$. Let $Tbra_{\langle n}$ denote the subset of $Tbra_n$ consisting of strings with $b_1=\langle$. We define {\it associative transformation for triple bra-ket} as follows:
\begin{df}
For any $b\in Tbra_{\langle n}$, suppose that $b_i=b_j=|$, and $b_k=\rangle$ with $i<j<k$. Let $b'=(b'_m)_{1\leq m\leq 2n+1}$ be the string of bracket defined by
\begin{displaymath}
b'_m=
\left\{
\begin{array}{ll}
b_{m+1}, & 1\leq m \leq i-2, \\
\langle, & m =i-1\\ 
b_{m+1}, & i\leq m \leq j-2,\\
|, & m=j-1,\\
b_{m+1}, & j\leq m\leq k-2,\\
|, & m=k-1,\\
b_{m+1}, & k\leq m\leq 2n,\\
\rangle, & m=2n+1.
\end{array}
\right.
\end{displaymath}
The string $b'$ is denoted by $ass_{\langle ||\rangle}(b)$, and is called the {\rm associative transformation for triple bra-ket} of $b$. 
\end{df}
\noindent
As is seen from the definition, this transformation is a natural generalization of that for bra-ket introduced in Definition 9.1. Furthermore we define an associative transformation to quadruple bra-ket of $b\in Tbra_{\langle n}$ as follows:

\begin{df}
For any $b\in Tbra_{\langle n}$, suppose that $b_i=|$, $b_j=|$, and $b_k=\rangle$. Let $b"=(b"_m)_{1\leq m\leq 2n+3}$ be the string of bracket defined by
\begin{displaymath}
b"_m=
\left\{
\begin{array}{ll}
\langle, & m=1,\\
b_m, & 2\leq m \leq i-1, \\
|, & m =i\\ 
b_m, & i+1\leq m \leq j-1,\\
|, & m=j,\\
b_m, & j+1\leq m\leq k-1,\\
|, & m=k,\\
b_m, & k+1\leq m\leq 2n+2,\\
\rangle, & m=2n+3.
\end{array}
\right.
\end{displaymath}
The string $b"$ is denoted by $ass_{to\langle |||\rangle}(b)$, and is called the {\rm associative transformation to quadruple bra-ket} of $b$.
\end{df}

\noindent
The following proposition provides us with a crucial ingredient to construct a family of linear subvarieties of $AC_N$ with $N\equiv 2\pmod 4$:

\begin{prp}
For any integer $n\geq 1$ and for any string $b\in Tbra_{\langle n}$, let $b'=ass_{\langle ||\rangle}(b)$ and $b"=ass_{to\langle |||\rangle}(b)$. Then we have
\begin{eqnarray}
V(f_{b"}(2))\subset V(f_b(2))\cap V(f_{b'}^{(+1)}(2)).
\end{eqnarray}
\end{prp}

\noindent
{\it Proof}. Let $b_i=|$, $b_j=|$, and $b_k=\rangle$ with $i<j<k$. Then it follows from the definition that for any $b_m\in L(b)-\{b_1,b_i,b_j\}$, we have 
\begin{eqnarray*}
f_{b_m}=f_{b'_{m-1}}^{(+1)}=f_{b"_m}.
\end{eqnarray*}
Therefore we have

\begin{eqnarray}
f_b(2)&=&\{f_{b_m};b_m\in L(b)-\{b_1,b_i,b_j\}\}\cup \{f_{b_i}(2), f_{b_j}(2)\},\\
f_{b'}^{(+1)}(2)&=&\{f_{b_m};b_m\in L(b)-\{b_1,b_i,b_j\}\}\cup \{f_{b'_{j-1}}^{(+1)}(2), f_{b'_{k-1}}^{(+1)}(2)\},\\ \nonumber
f_{b"}(2)&=&\{f_{b_m};b_m\in L(b)-\{b_1,b_i,b_j\}\}\\
&&\cup \{f_{b"_i}(2), f_{b"_j}(2), f_{b"_k}(2) \}.
\end{eqnarray}

Recall that the polynomials $f_{b_i}(2),\cdots,  f_{b"_k}(2)$ on the right hand sides of the above equalities are defined as
\begin{eqnarray*}
f_{b_i}(2)&=&\sum_{p\in cont(b_1)}x_p\sum_{q\in cont(b_i)}x_q-2,\\
f_{b_j}(2)&=&\sum_{q\in cont(b_i)}x_q\sum_{r\in cont(b_j)}x_r-2,\\
f_{b'_{j-1}}^{(+1)}(2)&=&\sum_{q\in cont(b'_{i-1})}x_{q+1}\sum_{r\in cont(b'_{j-1})}x_{r+1}-2\\
&=&\sum_{q\in cont(b_i)}x_q\sum_{r\in cont(b_j)}x_r-2,\\
f_{b'_{k-1}}^{(+1)}(2)&=&\sum_{r\in cont(b'_{j-1})}x_{r+1}\sum_{s\in cont(b'_{k-1})}x_{s+1}-2\\
&=&\sum_{r\in cont(b_j)}x_r\sum_{s\in cont(b'_{k-1})}x_{s+1}-2,\\
f_{b"_i}(2)&=&\sum_{p\in cont(b_1)}x_p\sum_{q\in cont(b_i)}x_q-2,\\
f_{b"_j}(2)&=&\sum_{q\in cont(b_i)}x_q\sum_{r\in cont(b_j)}x_r-2,\\
f_{b"_k}(2)&=&\sum_{r\in cont(b_j)}x_r\sum_{s\in cont(b'_{k-1})}x_{s+1}-2.\\
\end{eqnarray*}
Combining these with (11.2)-(11.4), we see that the inclusion (11.1) holds true. This completes the proof. \qed\\

\noindent
This proposition enables us to obtain a family of linear subvarieties of $AC_N$ for an arbitrary integer $N\equiv 2\pmod 4$:

\begin{thm}
For any integer $n\geq 1$ and for any string $b\in Tbra_{\langle n}$, we have
\begin{eqnarray}
V(f_{ass_{to\langle|||\rangle}(b)}(2))\subset V_{[1,2n+1]}\cap V_{[2,2n+2]}\cap V_{[1,2n+2]}((-1)^n).
\end{eqnarray}
In particular if $n$ is an even integer $2m$, say, then we have
\begin{eqnarray}
V(f_{ass_{to\langle|||\rangle}(b)}(2))\subset AC_{4m+2}.
\end{eqnarray}
\end{thm}

\noindent
{\it Proof}. It follows from Proposition 10.2 that
\begin{eqnarray}
V(f_b(2))&\subset& V_{[1,2n+1]},\\
V(f_{ass_{\langle|\rangle}(b)}^{(+1)}(2))&\subset& V_{[2,2n+2]}.
\end{eqnarray}
Furthermore it follows from Proposition 10.1 that
\begin{eqnarray}
V(f_{ass_{to\langle|||\rangle}(b)}(2))\subset V_{[1,2n+2]}((-1)^n).
\end{eqnarray}
Hence combining (11.6)-(11.8) with Proposition 11.1, we have the inclusion (11.5). This completes the proof. \qed\\\\

\section{Examples of AC-polygons}
The regular $n$-gons (including stars) are obviously examples of AC-polygons. These are obtained from our viewpoint as follows. Let $\Delta=\{(x, \cdots, x)\in\mathbb{A}^n;x\in \mathbb{R}\}$, the diagonal of $\mathbb{A}^n$. Then the intersection $V(u[1,n])\cap\Delta$ gives rise to the regular $n$-gons:

\begin{prp}
For any $n\geq 3$ and for any $k\in[0,n-1]$, let $c_{(n;k)}=\cos\frac{2k}{n}\pi$. Then we have
\begin{eqnarray*}
AC_n\cap\Delta=\{(c_{(n;k)},\cdots,c_{(n;k)});k\in [0,n-1]\}.
\end{eqnarray*}
\end{prp}

\noindent
{\it Proof}. As is explained in [, Proposition 2.7], we have 
\begin{eqnarray}
u[1,n]\Big|_{x_1=x,\cdots,x_n=x}=U_n(\frac{x}{2}),
\end{eqnarray}
where $U_n(z)$ denotes the Chebyshev polynomial of the second kind defined by
\begin{eqnarray}
U_n(z)=\frac{\sin (n+1)\theta}{\sin\theta}\hspace{2mm}\mbox{{\it with}}\hspace{2mm}z=\cos \theta.
\end{eqnarray}
Since $AC_n$ is defined to be $V(u[1,n]-1,u[1,n-1],u[2,n])$, it follows from (12.1) that
\begin{eqnarray*}
AC_n\cap\Delta=V(U_{n-1}\left(\frac{x}{2}\right), U_n\left(\frac{x}{2}-1\right)\}.
\end{eqnarray*}
Let us put $x=2\cos\theta$, then it follows from (12.2) that
\begin{eqnarray*}
U_{n-1}\left(\frac{x}{2}\right)&=&\frac{\sin n\theta}{\sin\theta},\\
U_n\left(\frac{x}{2}\right)&=&\frac{\sin (n+1)\theta}{\sin\theta}.
\end{eqnarray*}
Hence by a simple computation we see that the equality in Proposition 12.1 holds true. \qed\\

\noindent
Now for any integer $m$, we have by the addition formulas
\begin{eqnarray*}
&&2\cos\theta(\cos m\theta,\sin m\theta)-(\cos (m-1)\theta,\sin (m-1)\theta)\\
&=&(\cos (m+1)\theta,\sin (m+1)\theta).
\end{eqnarray*}
Hence if we take the initial points $p_0, p_1$ as
\begin{eqnarray*}
p_0&=&(1,0),\\
p_1&=&(\cos\frac{2k}{n}\pi,\sin\frac{2k}{n}\pi),
\end{eqnarray*}
and set $x_1=\cdots=x_n=c_{(n;k)}=2\cos\frac{2k}{n}\pi$, then we have
\begin{eqnarray*}
p_m&=&(\cos\frac{2km}{n}\pi,\sin\frac{2km}{n}\pi)
\end{eqnarray*}
for any $m\in\mathbb{Z}$. Thus the $n$-gon $P=(p_0,\cdots,p_{n-1})$ corresponding to the point $(c_{(n;k)},\cdots,c_{(n;k)})$ gives rise to the regular $n$-star $\{n/k\}$ in Schl\"{a}fli symbol.

\end{document}